\input amstex
\input amsppt.sty
\magnification=\magstep1
\hsize=160truemm
\vsize=229truemm
\baselineskip=16truept
\NoBlackBoxes
\nologo
\pageno=1
\topmatter
\TagsOnRight

\def\N{\Bbb N}
\def\Z{\Bbb Z}
\def\Q{\Bbb Q}

\def\l{\left}
\def\r{\right}
\def\b{\bigg}

\def\({\b(}
\def\[{\b[}
\def\){\b)}
\def\]{\b]}

\def\t{\text}
\def\f{\frac}

\def\ord{\roman{ord}}

\def\em{\emptyset}
\def\se {\subseteq}

\def\sm{\setminus}

\def\bi{\binom}
\def\eq{\equiv}

\def\ls{\leqslant}
\def\gs{\geqslant}

\def\ve{\varepsilon}
\def\da{\delta}

\def\Proof{\noindent{\it Proof}}
\def\Remark{\noindent{\it Remark}}

\def\Ack{\noindent {\bf Acknowledgments}}
\hbox{Sci. China Math. 63 (2020), no.\,3, 501--520.}
\medskip
\title Universal sums of three quadratic polynomials\endtitle
\author Zhi-Wei Sun \endauthor
\affil Department of Mathematics, Nanjing University
\\ Nanjing 210093, People's Republic of China
    \\  zwsun\@nju.edu.cn
    \\ {\tt http://maths.nju.edu.cn/$\sim$zwsun}
 \endaffil
\abstract  Let $a,b,c,d,e,f$ be integers with $a\gs c\gs e>0$, $b>-a$ and $b\eq a\pmod2$, $d>-c$ and $d\eq c\pmod 2$, $f>-e$ and $f\eq e\pmod2$.
Suppose that $b\gs d$ if $a=c$, and $d\gs f$ if $c=e$. When $b(a-b)$, $d(c-d)$ and $f(e-f)$ are not all zero, we prove that if each $n\in\N=\{0,1,2,\ldots\}$ can be written $x(ax+b)/2+y(cy+d)/2+z(ez+f)/2$
with $x,y,z\in\N$ then the tuple $(a,b,c,d,e,f)$ must be on our list of $473$ candidates, and show that 56 of them meet our purpose.
When $b\in[0,a)$, $d\in[0,c)$ and $f\in[0,e)$, we investigate the universal tuples $(a,b,c,d,e,f)$ over $\Z$ for which any $n\in\N$ can be written $x(ax+b)/2+y(cy+d)/2+z(ez+f)/2$
with $x,y,z\in\Z$, and show that there are totally 12082 such candidates some of which are proved to be universal tuples over $\Z$.
For example, we show that any $n\in\N$ can be written as $x(x+1)/2+y(3y+1)/2+z(5z+1)/2$ with $x,y,z\in\Z$, and conjecture that each $n\in\N$ can be written as $x(x+1)/2+y(3y+1)/2+z(5z+1)/2$ with $x,y,z\in\N$.

\endabstract
\thanks 2010 {\it Mathematics Subject Classification}.
Primary 11E25; Secondary 11D85, 11E20.
\newline \indent Supported by the National Natural Science
Foundation of China (grant 11571162) and the NSFC-RFBR Cooperation and Exchange Program (grant 11811530072).
\endthanks
\endtopmatter
\document

\heading{1. Introduction}\endheading

Recall that the triangular numbers have the form $T_x=x(x+1)/2$ with $x\in\Z$. Since $T_{-1-n}=T_n$ for all $n\in\N=\{0,1,2,\ldots\}$, we have $\{T_x:\ x\in\Z\}=\{T_n:\ n\in\N\}$.
In 1796 Gauss proved Fermat's assertion that each $n\in\N$ can be written as the sum of three triangular numbers.
For each $m=3,4,5,\ldots$ the {\it $m$-gonal numbers} (or polygonal numbers of order $m$) are given by
$$p_m(x):=(m-2)\bi x2 +x=\f{x((m-2)x-(m-4))}2\ \ (x\in\N).$$
Note that $p_3(x)=T_x$ and $p_4(x)=x^2$. For $m\in\{5,6,7,\ldots\}$, those $p_m(x)$ with $x\in\Z$ are called {\it generalized $m$-gonal numbers}.

For a subset $S$ of $\Z$ and polynomials $f_1(x),f_2(x),f_3(x)$ with $f_i(S)=\{f_i(x):\ x\in S\}\se\N$ for $i=1,2,3$, if
any $n\in\N$ can be written as $f_1(x)+f_2(y)+f_3(z)$ with $x,y,z\in S$ then we call the sum
$f_1(x)+f_2(y)+f_3(z)$ {\it universal over $S$}.

Let $\Z^+=\{1,2,3,\ldots\}$. In 1862 Liouville (cf. Berndt [B, p.\,82] and Dickson [D99, p.\,23])
determined all those universal sums $aT_x+bT_y+cT_z$ (over $\N$ or $\Z$) with $a,b,c\in\Z^+$.
It is known that $ax^2+by^2+cz^2$ is not universal (over $\N$ or $\Z$) for any $a,b,c\in\Z^+$ (cf. [DW]).
The determination of those universal sums $ap_i(x)+bp_j(y)+cp_k(z)$ (over $\N$ or $\Z$) with $\{i,j,k\}=\{3,4\}$ and $a,b,c\in\Z^+$
was proposed by the author [S07] and completed via the three papers [S07], [GPS] and [OS]. Note that
$$\{T_x+T_y:\ x,y\in\Z\}=\{x^2+2T_y:\ x,y\in\Z\}\tag1.1$$
as observed by Euler (cf. [D99, p.\,11]); in fact, $x^2+2T_y=T_{x+y}+T_{y-x}=T_{x+y}+T_{x-y-1}.$

The author [S15] showed that there are only 95 candidates for universal sums $ap_i(x)+bp_j(y)+cp_k(z)$ over $\N$ with $a,b,c\in\Z^+$, $i,j,k\in\{3,4,5,\ldots\}$ and $\max\{i,j,k\}\gs5$.
Though none of the 95 sums has been proved to be universal over $\N$, many of them have been proved to be universal over $\Z$ (cf. [S15] and [JOS]).

For $c\in\Z^+$ and $m\in\{3,4,\ldots\}$, clearly
$cp_m(x)=x(a_0x+b_0)/2$ with $a_0=c(m-2)$ and $b_0=-c(m-4)\in(-a_0,a_0]$.
Instead of $cp_m(x)$, we may consider more general polynomials
$$\psi_{a,b}(x):=\f{x(ax+b)}2\ \t{with}\ a\in\Z^+,\ b\in\Z,\ b>-a\ \t{and}\ a\eq b\pmod 2.\tag1.2$$
Clearly, $\psi_{a,b}(\N)\se\N.$
For positive integers $a,c,e$ and integers $b>-a,\ d>-c,\ f>-e$ with $a+b,c+d,e+f$ all even, if $\psi_{a,b}(x)+\psi_{c,d}(y)+\psi_{e,f}(z)$
is universal over $\N$ then we simply call the ordered tuple $(a,b,c,d,e,f)$ {\it universal over $\N$}.
In view of Liouville's result (cf. [D99, p.\,23]), and [S07], [GPS] and [OS], all those universal tuples
$(a,b,c,d,e,f)$ over $\N$ with $b\in\{0,a\}$, $d\in\{0,c\}$ and $f\in\{0,e\}$
have been determined.

In our first theorem, we give some new universal tuples $(a,b,c,d,e,f)$ over $\N$ with $a\mid b$, $c\mid d$ and $e\mid f$.
\proclaim{Theorem 1.1} All the following $56$ ordered tuples
$$\align &(1, 3, 1, 1, 1, 1),\ (1, 3, 1, 3, 1, 1),\ (1, 5, 1, 1, 1, 1),\ (1, 5, 1, 3, 1, 1),\ (1, 7, 1, 1, 1, 1),
\\& (1, 7, 1, 3, 1, 1),\ (1, 9, 1, 1, 1, 1),\ (2, 0, 1, 3, 1, 1),\ (2, 0, 1, 3, 1, 3),\ (2, 0, 1, 5, 1, 1),
\\&(2, 0, 1, 5, 1, 3),\ (2, 0, 1, 7, 1, 1),\ (2, 0, 1, 7, 1, 3),\ (2, 0, 1, 9, 1, 1),\ (2, 0, 1, 9, 1, 3),
\\& (2, 0, 1, 11, 1, 1),\ (2, 0, 1, 11, 1, 3),\ (2, 0, 1, 13, 1, 1),\ (2, 0, 1, 13, 1, 3),\ (2, 0, 1, 15, 1, 1),
\\& (2, 0, 2, 0, 1, 3),\ (2, 2, 1, 3, 1, 1),\ (2, 2, 1, 5, 1, 1),\ (2, 2, 1, 7, 1, 1),\ (2, 2, 2, 0, 1, 3),
\\& (2, 2, 2, 0, 1, 5),\ (2, 2, 2, 0, 1, 7),\ (2, 2, 2, 0, 1, 9),\ (2, 4, 1, 1, 1, 1),\ (2, 4, 2, 0, 1, 1),
\\& (2, 4, 2, 0, 1, 3),\ (2, 4, 2, 2, 1, 1),\ (2, 4, 2, 2, 2, 0),\ (2, 6, 1, 1, 1, 1),\ (2, 6, 1, 3, 1, 1),
\\& (2, 6, 2, 0, 1, 1),\ (2, 6, 2, 0, 1, 3),\ (2, 6, 2, 2, 1, 1),\ (2, 6, 2, 2, 2, 0),\ (2, 8, 1, 1, 1, 1),
\\&(2, 8, 2, 0, 1, 1),\ (2, 8, 2, 0, 1, 3),\ (2, 8, 2, 2, 2, 0),\ (2, 10, 2, 0, 1, 1),\ (2, 10, 2, 0, 1, 3),
\\& (2, 12, 2, 0, 1, 1),\ (2, 12, 2, 0, 1, 3),\ (2, 14, 2, 0, 1, 1),\ (3,3,2,0,1,3),\ (3,9,2,0,1,1),
\\& (3,9,2,0,1,3),\ (4,0,1,3,1,1),\ (4,0,1,5,1,1),\ (4,0,1,7,1,1),\ (4,4,1,3,1,1)
\\& (8,0,1,3,1,1)
\endalign$$
are universal over $\N$.
\endproclaim

 We have the following conjecture on other possible universal tuples $(a,b,c,d,e,f)$ over $\N$ with $a\mid b$, $c\mid d$ and $e\mid f$.

\proclaim{Conjecture 1.2} The following $10$ ordered tuples
$$\align &(4,0,2,0,1,3),\ (4,0,2,0,1,5),\ (4,0,2,6,1,1),\ (4,0,2,6,2,0),\ (4,4,2,0,1,3),
\\& (4,8,2,0,1,1),\ (4,8,2,0,1,3),\ (4,12,2,0,1,1),\ (6,0,2,0,1,3),\ (6,6,2,0,1,3)
\endalign$$
are universal over $\N$.
\endproclaim
\Remark\ 1.3. It is easy to see that $(4,8,2,0,1,3)$ is universal over $\N$ if and only if any integer $n>2$ can be written as $x^2+2y^2+T_z$ with $x\in\N$ and $y,z\in\Z^+$.
\medskip

 Now we state our second theorem.

\proclaim{Theorem 1.4} Let $a,c,e$ be positive integers and let
$b>-a$, $d>-c$ and $f>-e$ be integers with $a+b,c+d,e+f$ all even.
Suppose that $a\gs c\gs e$, and $b\gs d$ if $a=c$, and $d\gs f$ if $c=e$, and that the ordered tuple $(a,b,c,d,e,f)$ is universal over $\N$.

{\rm (i)} If $a\mid b$, $c\mid d$ and $e\mid f$, but $b(a-b),d(c-d),f(e-f)$ are not all zero, then $(a,b,c,d,e,f)$ must be among the $56$ tuples in Theorem 1.1 or the $10$ tuples in Conjecture 1.1.

{\rm (ii)} If $a\nmid b$ or $c\nmid d$ or $e\nmid f$, then $(a,b,c,d,e,f)$ must be among the $407$ tuples listed in the Appendix.
\endproclaim

\proclaim{Conjecture 1.5} All the $407$ tuples in the Appendix are universal over $\N$. In particular,
$$\l\{\f{x(x+1)}2+\f{y(3y+1)}2+\f{z(5z+1)}2:\ x,y,z\in\N\r\}=\N.\tag1.3$$
\endproclaim
\Remark\ 1.6. The author would like to offer 135 US dollars for the first proof of (1.3).
In [S17] the author conjectured that any $n\in\Z^+$ can be written as $x^3+y^2+T_z$ with $x,y\in\N$ and $z\in\Z^+$;
we also conjecture that $x^3+T_y+z(3z+7)/2$ is universal over $\N$.
\medskip

Guy [G94] noted that $p_5(x)+p_5(y)+p_5(z)$ is universal over $\Z$.
Sun [S15] proved that if $ap_5(x)+bp_5(y)+cp_5(z)$ is universal over $\Z$ with $a,b,c\in\Z^+$ and $a\ls b\ls c$
then $(a,b,c)$ is among the following 20 triples:
$$(1,1,i)\ (i=1,\ldots,6,8,9,10),\ (1,2,j)\ (j=2,3,4,6,8),\ (1,3,k)\ (k=3,4,6,7,8,9).$$
By Sun [S15, Theorem 1.1(ii)], Ge and Sun [GS], and Oh [O11], for each $(a,b,c)$ among the 20 triples the sum $ap_5(x)+bp_5(y)+cp_5(z)$ is indeed universal over $\Z$.

Recall that $T_x+T_y+T_z$ is universal over $\Z$. It is easy to see that
$$\psi_{1,1}(\Z)=\l\{T_x:\ x\in\Z\r\}=\{x(2x+1):\ x\in\Z\}=\psi_{4,2}(\Z).\tag1.4$$
Motivated by this, Sun [S17] investigated universal sums over $\Z$
of the form $x(ax+1)+y(by+1)+z(cz+1)$ with $1\ls a\ls b\ls c$, or the form $x(ax+b)+y(ay+c)+z(az+d)$ with $2<a\gs b\gs c\gs d\gs0$.
Later, Ju and Oh [JO] proved some conjectures of Sun [S17] in this direction.

For $a\in\Z^+$, clearly $\psi_{a,-b}(\Z)=\psi_{a,b}(\Z)$ for all $b=0,\ldots,a$ with $b\eq a\pmod2$, and $\psi_{a,a}(\Z)=\psi_{4a,2a}(\Z)$ by (1.4).
Thus we are led to find all the sums
$$\psi_{a,b}(x)+\psi_{c,d}(y)+\psi_{e,f}(z)=\f{x(ax+b)}2+\f{y(cy+d)}2+\f{z(ez+f)}2\tag1.5$$
which are universal over $\Z$, where $a,c,e\in\Z^+$, $b,d,f\in\N$, $b<a$ and $a\eq b\pmod2$, $d<c$ and $c\eq d\pmod2$, and $f<e$ and $e\eq f\pmod2$. If the sum in (1.5) is universal over $\Z$, then we say that the ordered tuple $(a,b,c,d,e,f)$ is universal over $\Z$.

\proclaim{Theorem 1.7} Let $a,b,c,d,e,f\in\N$ with $a>b$, $c>d$, $e>f$, $a\eq b\pmod2$, $c\eq d\pmod 2$, $e\eq f\pmod2$, $a\gs c\gs e\gs2$,
and $b\gs d$ if $a=c$, and $d\gs f$ if $c=e$. Suppose that the ordered tuple $(a,b,c,d,e,f)$ is universal over $\Z$. Then $(a,b,c,d,e,f)$ must be among the $12082$ ordered tuples listed in [S17a].
\endproclaim
\Remark\ 1.8. Chan and Oh [CO] showed that there are only finitely many
equivalence classes of positive ternary universal integral quadratic polynomials. We have analysed those tuples $(a,b,c,d,e,f)$ with $a\ls5$ listed in [S17a], only the following 10 tuples
$$\align&(5,1,2,0,2,0),\ (5,3,2,0,2,0),\ (5,1,4,0,2,0),\ (5,3,4,0,2,0),\ (5,1,4,0,3,1),
\\&(5,3,4,0,3,1),\ (5,1,5,1,2,0),\ (5,3,5,3,2,0),\ (5,3,4,0,4,0),\ (5,3,5,3,4,0)
\endalign$$
have not yet been proved to be universal over $\Z$.
\medskip

For polynomials $f_1,f_2,f_3,f_4$ with $f_i(\Z)\se\N$, if
$$\{f_1(x)+f_2(y):\ x,y\in\Z\}=\{f_3(x)+f_4(y):\ x,y\in\Z\}$$
then we say that $f_1(x)+f_2(y)$ is equivalent to $f_3(x)+f_4(y)$ and write $f_1(x)+f_2(y)\sim f_3(x)+f_4(y)$ for this.
(1.1) indicates that $T_x+T_y\sim x^2+2T_y$.
In light of this, we obtain the following auxiliary result which has its own interest.

\proclaim{Theorem 1.9} {\rm (i)} For any $a\in\Z^+$ and $b\in\N$ with $b\ls a/2$, we have
$$x(ax+b)+y(ay+a-b)\sim aT_x+ \psi_{a,a-2b}(y).\tag1.6$$

{\rm (ii)} We have
$$\align x^2+T_y&\sim p_5(x)+2p_5(y),\tag1.7
\\T_x+2T_y&\sim p_5(x)+p_8(y),\tag1.8
\\x^2+4T_y&\sim 4p_5(x)+p_8(y),\tag1.9
\\T_x+T_y&\sim \psi_{5,1}(x)+\psi_{5,3}(y).\tag1.10
\endalign$$
\endproclaim
\Remark\ 1.10. (1.6) with $a=1$ and $b=0$ gives (1.1). Putting $a=3$ and $b=1$ in (1.6) we get
$$2p_5(x)+p_8(y)\sim 3T_x+p_5(y).\tag1.11$$

With helps of Theorem 1.9 and the theory of ternary quadratic forms, we establish the following new result.

\proclaim{Theorem 1.11} {\rm (i)} We have
$$\l\{\f{x(x+1)}2+\f{y(3y+1)}2+\f{z(5z+1)}2:\ x,y,z\in\Z\r\}=\N.\tag1.12$$

{\rm (ii)} For any $\da\in\{0,1\}$ and $r\in\{1,3,5\}$, we have
$$\l\{x(x+\da)+\f{y(3y+1)}2+\f{z(5z+r)}2:\ x,y,z\in\Z\r\}=\N.\tag1.13$$

{\rm (iii)} For any  $r,s,t\in\{1,3\}$ with $\{r,s\}\not=\{3\}$, we have
$$\l\{\f{x(3x+r)}2+\f{y(3y+s)}2+\f{z(5z+t)}2:\ x,y,z\in\Z\r\}=\N.\tag1.14$$

{\rm (iv)} Let $s,t\in\{1,3,5\}$ with $\{s,t\}\not=\{5\}$. Then
$$\l\{\f{x(3x+1)}2+\f{y(5y+s)}2+\f{z(5z+t)}2:\ x,y,z\in\Z\r\}=\N.\tag1.15$$
\endproclaim
\Remark\ 1.12. The author [S17, Conjecture 1.2] conjectured that $x^2+y(3y+1)/2+z(5z+3)/2$ is universal over $\N$.
In view of (1.11), and (1.14) with $r=3$, $s=1$ and $t\in\{1,3\}$, the tuples $(6,4,6,2,5,1)$ and $(6,4,6,2,5,3)$ are universal over $\Z$.
\medskip

We are going to prove Theorem 1.1 in Section 3 based on some lemmas given in the next section.
Theorems 1.4 and 1.7 and Theorems 1.9 and 1.11 will be shown in Sections 4 and 5 respectively.

\heading{2. Some lemmas}\endheading

The Gauss-Legendre theorem on sums of three squares (cf. [N96, pp.\,17-23]) asserts that $\{x^2+y^2+z^2:\,x,y,z\in\Z\}=\N\sm\{4^k(8l+7):\,k,l\in\N\}$.
For $n\in\N$ we define
$$R_3(n):=|\{(x,y,z)\in\Z^3:\ x^2+y^2+z^2=n\ \t{and}\ \gcd(x,y,z)=1\}|.$$

\proclaim{Lemma 2.1 {\rm (Gauss)}} Let $n\in\N$. Then $R_3(1)=6,\ R_3(2)=12$, $R_3(3)=8$, and
$$R_3(n)=\cases 12h(-n)&\t{if}\ n>3\ \t{and}\ n\eq1,2\pmod4,
\\24h(-n)&\t{if}\ n>3\ \t{and}\ n\eq3\pmod8,
\\0&\t{if}\ 4\mid n\ \t{or}\ n\eq7\pmod 8,
\endcases$$
where $h(-n)$ denotes the class number of the field $\Q(\sqrt{-n})$.
\endproclaim
\Remark\ 2.2. One may consult [P, p.\,140] for this classical result.

\proclaim{Lemma 2.3} Let $n\in\Z^+$ be squarefree.

{\rm (i)} $h(-n)=1$ if and only if $n$ is among the following nine numbers
$$1,\ 2,\ 3,\ 7,\ 11,\ 19,\ 43,\ 67,\ 163.$$

{\rm (ii)} $h(-n)=2$ if and only if $n$ is among the following numbers
$$5,\ 6,\ 10,\ 13,\ 15,\ 22,\ 35,\ 37,\ 51,\ 58,\ 91,\ 115,\ 123,\ 187,\ 235,\ 267,\ 403,\ 427.$$

{\rm (iii)} $h(-n)=3$ if and only if $n$ is among the following numbers
$$23,\ 31,\ 59,\ 83,\ 107,\ 139,\ 211,\ 283,\ 307,\ 331,\ 379,\ 499,\ 547,\ 643,\ 883,\ 907.$$

{\rm (iv)} $h(-n)=4$ if and only if $n$ is among the following numbers
$$\align &14,\ 17,\ 21,\ 30,\ 33,\ 34,\ 39,\ 42,\ 46,\ 55,\ 57,\ 70,\ 73,\ 78,\ 82,\ 85,\ 93,\ 97,
\\&102,\ 130,\ 133,\ 142,\ 155,\ 177,\ 190,\ 193,\ 195,\ 203,\ 219,\ 253,\ 259,\ 291,
\\&323,\ 355,\ 435,\ 483,\ 555,\ 595,\ 627,\ 667,\ 715,\ 723,\ 763,\ 795,\ 955,\ 1003,
\\&1027,\ 1227,\ 1243,\ 1387,\ 1411,\ 1435,\ 1507,\ 1555.
\endalign$$

{\rm (v)} If $n\eq 1\pmod 8$, then $h(-n)\in\{5,6,7,8\}$ if and only if $n$ is among the following numbers
$$41,\,65,\, 105,\, 113,\, 137,\, 145,\, 217,\, 265,\, 273,
\, 313,\, 337,\, 345,\, 385,\, 457,\, 505,\, 553,\, 697,\, 793.$$
\endproclaim
\Remark\ 2.4. This is a known result, see,  [A], [ARW], [W], [Wa] and [We].

\proclaim{Lemma 2.5} Any integer $n>1$ can be written as $T_x+T_y+T_z$ with $x,y\in\Z^+$ and $z\in\N$.
\endproclaim
\Proof. By Gauss' result, $n=T_x+T_y+T_z$ for some $x,y,z\in\N$ with $x\gs y\gs z$. If $y>0$ then $x,y\in\Z^+$.
If $y=0$, then $z$ is also zero and $n=T_x$ is a triangular number.

Now assume that $n=T_m>1$ with $m\in\N$. Clearly, $m>1$. If $8T_m+3=(2m+1)^2+2$  is divisible by $d^2$ for some integer $d>1$, then
$$|\{(x,y,z)\in\Z^3:\ x^2+y^2+z^2=8T_m+3\ \&\ \gcd(x,y,z)=d\}|
=R_3\l(\f{8T_m+3}{d^2}\r)>0$$
by Lemma 2.1, and hence $8T_m+3=x^2+y^2+z^2$ for some positive odd integers $x,y,z$ with $\min\{x,y\}\gs d>1$.
When $8T_m+3=(2m+1)^2+2$ is squarefree, by Lemma 2.1 and Lemma 2.3(i) we have
$R_3(8T_m+3)=24h(-8T_m-3)>24$
since $8T_m+3=(2m+1)^2+2\not\in\{1,2,3,7,11,19,43,67,163\}$, hence there are integral solutions of the equation
$x^2+y^2+z^2=8T_m+3$ other than the 24 trivial solutions
$$(\pm1,\pm1,\pm(2m+1)),\ (\pm1,\pm(2m+1),\pm1),\ (\pm(2m+1),\pm1,\pm1).$$
Thus, for some $x,y\in\Z^+$ and $z\in\N$ we have
$$8T_m+3=(2x+1)^2+(2y+1)^2+(2z+1)^2, \ \t{i.e.},\ T_m=T_x+T_y+T_z.$$
This concludes our proof. \qed

\proclaim{Lemma 2.6} Let $n>2$ be an integer. Then
$$4n+2=x^2+y^2+z^2 \ \ \ \t{for some}\ x,y,z\in\N\ \t{with}\ x\gs y>1.$$
\endproclaim
\Proof. If $d^2\mid 4n+2$ for some integer $d>1$, then $R_3((4n+2)/d^2)>0$ (by Lemma 2.1) and hence $4n+2=(dx)^2+(dy)^2+(dz)^2$
for some $x,y,z\in\N$ with two of $dx,dy,dz$ odd and greater than one.
If neither $4n$ nor $4n+1$ is a square, then $4n+2$ is not of the form $x^2+2$ or $x^2+1$ with $x\in\Z$,
and by the Gauss-Legendre theorem we can write $4n+2=x^2+y^2+z^2$ with $x,y,z\in\N$ and $x\gs y>1$.

Now assume that $4n+\da=t^2$ with $\da\in\{0,1\}$ and $t\in\N$, and that $4n+2$ is squarefree. Then
$R_3(4n+2)=12h(-4n-2)>24$ by Lemmas 2.1 and 2.2. (Note that none of $22-2,\,22-1,\,58-2,\,58-1$ is a square.)
So the equation
$4n+2=x^2+y^2+z^2$
has integral solutions other than the trivial solutions
$$(\pm1,\pm\bar\da,\pm t),\ (\pm\bar\da,\pm1,\pm t),(\pm1,\pm t,\pm\bar\da),(\pm\bar\da,\pm t,\pm1),\ (\pm t,\pm1,\pm\bar\da),\ \ (\pm t,\pm\bar\da,\pm1),$$
where $\bar\da=1-\da$. (No matter $\da=0$ or $1$ there are exactly 24 trivial solutions.)
Hence $4n+2=x^2+y^2+z^2$ for some $x,y,z\in\N$ with $x\gs y\gs z$ and $y>1$.
This ends the proof. \qed

\proclaim{Lemma 2.7} Let $n>1$ be an odd integer. Then, for each $m=1,2$ there are $x,y,z\in\N$ with $\max\{x,z\}>0$ and $\max\{y,z\}>0$ such that $x^2+y^2+mz^2=n^2$.
\endproclaim
\Proof. In 1907 Hurwitz (cf. [D99, p.\,271]) showed that
$$\align |\{(x,y,z)\in\Z^3:\ x^2+y^2+z^2=n^2\}|
=6\prod_{p\mid n}\(p^{\ord_p(n)}+\l(1-\l(\f{-1}p\r)\r)\f{p^{\ord_p(n)}-1}{p-1}\),
\endalign$$
where $\ord_p(n)$ is the order of $n$ at the prime $p$ and $(\f{\cdot}p)$ is the Legendre symbol.
This implies that
$$|\{(x,y,z)\in\Z^3:\ x^2+y^2+z^2=n^2\}|\gs 6\prod_{p\mid n}p^{\ord_p(n)}=6n>6.$$
So the equation $x^2+y^2+z^2=n^2$ has integral solutions other than the 6 trivial solutions
$(0,0,\pm n),\ (0,\pm n,0),\ (\pm n,0,0).$
This proves the desired result in the case $m=1$.

Now we consider the case $m=2$. By a result of Cooper and Lam [CL], we have
$$\align |\{(x,y,z)\in\Z^3:\ x^2+y^2+2z^2=n^2\}|
=&4\prod_{p\mid n}\f{p^{\ord_p(n)+1}-1-(\f{-2}p)(p^{\ord_p(n)}-1)}{p-1}
\\\gs&4\prod_{p\mid n}\f{p^{\ord_p(n)+1}-p^{\ord_p(n)}}{p-1}=4n>4.
\endalign$$
So the equation $x^2+y^2+2z^2=n^2$ has integral solutions other than the 4 trivial solutions
$(\pm n,0,0)$ and $(0,\pm n,0)$.
This proves the desired result in the case $m=2$. \qed

\proclaim{Lemma 2.8} Let $m>3$ be an integer. Then we can write $T_m=x^2+y^2+T_z$ with $x,y,z\in\N$, $y\gs3$ and $z\gs1$.
\endproclaim
\Proof. It is easy to verify the desired result for $m=4,5,6,7,8$.
When $m=a^2$ with $a\in\{3,4,\ldots\}$, the desired result also holds since
$T_m=m+T_{m-1}=0^2+a^2+T_{a^2-1}.$

Now we assume that $m$ is greater than 8 and not a square.
By Dickson [D39, pp.\,112-113],
$$\N\sm\{x^2+y^2+2z^2:\ x,y,z\in\Z\}=\{4^k(16l+14):\ k,l\in\N\}.$$
Thus, in view of Lemma 2.7, no matter $2m+1$ is a square or not,  $2m+1=x^2+y^2+2z^2$ for some $x,y,z\in\N$ with $(x^2+z^2)(y^2+z^2)>0$.
Note that $x^2+y^2>1$.
If $(x^2+z^2)(y^2+z^2)\ls 16$, then $xy\ls 4$ and $z\ls 2$, hence $2m+1=x^2+y^2+2z^2\ls 17$
which contradicts $m>8$. Thus $(x^2+z^2)(y^2+z^2)>16$.

Observe that
$$\align 8T_m+1=&(2m+1)^2=(x^2+z^2+(y^2+z^2))^2
\\=&(x^2+z^2-(y^2+z^2))^2+4(x^2+z^2)(y^2+z^2)
\\=&(x^2-y^2)^2+4((xy+z^2)^2+(xz-yz)^2).
\endalign$$
Since $x\not\eq y\pmod2$ and $|x^2-y^2|\gs x+y>1$, we have $|x^2-y^2|=2w+1$ for some $w\in\Z^+$.
Note that $(xy+z^2)^2+(|x-y|z)^2=(x^2+z^2)(y^2+z^2)\eq0\pmod2$. Thus
$$u:=\f{xy+z^2+|x-y|z}2\in\N\ \ \t{and}\ \ v:=\l|\f{xy+z^2-|x-y|z}2\r|\in\N.$$
Since
$$u^2+v^2=\f{(xy+z^2)^2+(xz-yz)^2}2=\f{(x^2+z^2)(y^2+z^2)}2>8,$$
we have $\max\{u,v\}\gs3$. Finally,
$8T_m+1=(2w+1)^2+8(u^2+v^2)$ and hence $T_m=u^2+v^2+T_w.$
This concludes the proof. \qed

\heading{3. Proof of Theorem 1.1}\endheading

\medskip
\noindent{\it Proof of Theorem 1.1}. We first make a useful observation:
$$\f {z(z+2k+1)}2=T_{z+k}-T_k\ \ \t{for all}\ k\in\N.$$

(a) Let $k\in\{1,2,3\}$. Clearly $(1,2k+1,1,3,1,1)$ is universal over $\N$ if and only if
each integer $n\gs T_k+1$ can be written as $T_{x+k}+T_{y+1}+T_z$ with $x,y,z\in\N$.
Let $n$ be any integer with $n\gs T_k+1$.
By Lemma 2.3, we can write $n$ as $T_x+T_y+T_z$ with $x,y,z\in\N$ and $x\gs y\gs \max\{z,1\}$.
If $x<k$, then $k\in\{2,3\}$ and $n=T_x+T_y+T_z\ls 3T_{k-1}$. Note that $3T_{2-1}<T_2+1$, $T_3+1=7$ and $3T_{3-1}=9$.
Clearly, $7=T_3+T_1+T_0$, $8=T_3+T_1+T_1$ and $9=T_3+T_2+T_0.$
Thus $(1,2k+1,1,3,1,1)$ is universal over $\N$.

Similarly, for each $k\in\{1,2,3,4\}$, the tuple $(1,2k+1,1,1,1,1)$ is universal over $\N$ and so is $(2,2,2,0,1,2k+1)$
in view of $(1.1)$.

(b) By Lemma 2.6, for any integer $n>2$ there are $x,y\in\Z^+$ and $z\in\N$ with $z\ls y$ such that
$$4n+2=(2x+1)^2+(y+z+1)^2+(y-z)^2, \ \t{i.e.,}\ n=2T_x+T_y+T_z.$$
So, for any $n\in\N$ there are $x,y,z\in\N$ such that $n+3=2T_{x+1}+T_{y+1}+T_z$ and hence
$n=x(x+3)+y(y+3)/2+z(z+1)/2$. This prove the universality of $(2,6,1,3,1,1)$ over $\N$.
As $y(y+3)/2=T_{y+1}-1$ and $T_y+T_z\sim y^2+2T_z$, the tuples $(2,6,1,1,1,1)$ and $(2,6,2,2,2,0)$ are also universal over $\N$.

Let $k\in\{1,2,3\}$. It is easy to check that each $n=T_k,T_k+1,\ldots,k^2-k+6$ can be written as $2T_x+T_y+T_z$ with $x,y,z\in\N$ and $z\gs k$. Let $n\in\N$ with $n>k(k-1)+6$. Then $n=2T_x+T_y+T_z$ for some $x,y,z\in\N$. As $n>2T_2+T_{k-1}+T_{k-1}$,
either $x\gs 3$ or $\max\{y,z\}\gs k$.
If $\max\{y,z\}<k$, then $x\gs 3\gs k$,  $T_y+T_z\ls 2T_{k-1}=k(k-1)\ls 6$ and $T_y+T_z\not=5$.
It is easy to check that $2T_x+2\da=2T_{\da}+T_x+T_x$ and $2T_x+2\da+1=2T_{\da}+T_{x-1}+T_{x+1}$ for $\da=0,1$;
also, $2T_x+4=2T_0+T_{x-2}+T_{x+2}$ and $2T_x+6=T_0+T_3+2T_x$.
Therefore $(2,2,1,2k+1,1,1)$ is a universal tuple over $\N$.

In view of Lemma 2.6, for any $n\in\Z^+$ we can write
$$8n+6=4(2n+1)+2=(2x+1)^2+(2y+1)^2+w^2$$
with $x\in\Z^+$ and $y,w\in\N$. Since $w^2\eq 4\pmod 8$, we have $w=2(2z+1)$ for some $z\in\N$.
Therefore $8n+6=(2x+1)^2+(2y+1)^2+4(2z+1)^2$, hence $n=T_x+T_y+4T_z$ and $n-1=m(m+3)/2+T_y+4T_z$ with $m=x-1\in\N$.
This proves the universality of $(4,4,1,3,1,1)$ over $\N$.

(c) The tuple $(2,0,1,15,1,1)$ is universal over $\N$ if and only if any integer $n\gs 28$ can be written
in the form $x^2+T_y+T_z$ with $x,y,z\in\N$ and $\max\{y,z\}\gs7$. It is easy to verify that every $n=28,29,\ldots,78$ can be written as $x^2+T_y+T_z$ with $x,y,z\in\N$ and $z\gs7$.
Now let $n\in\N$ with $n>78$. We can write $n$ as $x^2+T_y+T_z$ with $x,y,z\in\N$ (cf. [S07]).
Suppose that $\max\{y,z\}\ls6$. Then
$n-x^2=T_y+T_z$ belongs to the set
$$R:=\{T_i+T_j:\, i,j=0,\ldots,6\}=(\{0,\ldots,31\}\cup\{36,42\})\setminus\{5,8,14,17,19,23,26,28,29\}$$
and hence $x>6$ since $36+42<n$. Note that
$$R\se\bigcup_{m=0}^6\{m^2,\,m^2+2,\,m(m+1)+0^2,\,\ldots,\,m(m+1)+5^2\}.$$
Also, $x^2=T_{x-1}+T_x$, $x^2+2= T_{x-2}+T_{x+1}$ and $x^2+m(m+1)=T_{x-m-1}+T_{x+m}$.
So the tuple $(2,0,1,15,1,1)$ is indeed universal over $\N$.

Let $k\in\{2,3,4,5,6\}$. By the last paragraph, any integer $n\gs T_7=28$ can be written in the form $x^2+T_y+T_z$
with $x,y,z\in\N$ and  $z\gs 7\gs k$. It is easy to verify that each $n=T_k,T_{k}+1,\ldots,27$ can be written as $x^2+T_y+T_z$ with $x,y,z\in\N$ and $z\gs k$.
Therefore the tuple $(2,0,1,2k+1,1,1)$ is universal over $\N$.

(d) Let $n>1$ be an integer. We claim that the equation $4n+1=x^2+y^2+z^2$ has at most $2^3\times3!=48$ integral solutions with $x-z,y-z\in\{\pm1\}$.

For any $y,z\in\Z^+$, it is easy to see that $2(y-1)^2+y^2\not=2(z+1)^2+z^2$
 by considering the cases $y<z$, $y\in\{z,z+1\}$, and $y>z+1$.

 Suppose that $4n+1=(z_0+1)^2+(z_0-1)^2+z_0^2$ for some $z_0\in\N$. Then $z_0>0$. For any $z\in\N$ with $z<z_0$ we have
 $$\align 2(z-1)^2+z^2<&(z+1)^2+(z-1)^2+z^2<2(z+1)^2+z^2\ls 2z_0^2+(z_0-1)^2
 \\<&(z_0+1)^2+(z_0-1)^2+z_0^2=4n+1,
 \endalign$$
 and  for any integer $z>z_0+1$ we have
 $$\align 2(z+1)^2+z^2\gs&(z+1)^2+(z-1)^2+z^2\gs 2(z-1)^2+z^2\gs 2(z_0+1)^2+(z_0+2)^2
 \\>&(z_0+1)^2+(z_0-1)^2+z_0^2=4n+1.
 \endalign$$
 Note also that
 $$\align 2(z_0-1)^2+z_0^2<&(z_0+1)^2+(z_0-1)^2+z_0^2=4n+1
 \\<&2z_0^2+(z_0+1)^2<2(z_0+1)^2+z_0^2<2(z_0+2)^2+(z_0+1)^2.
 \endalign$$
 So $4n+1\not=2(z\pm1)^2+z^2$ for any $z\in\N$.

In view of the above, the claim does hold.

  If $d^2\mid 4n+1$ for some integer $d>1$, then $R_3((4n+1)/d^2)>0$ by Lemma 2.1, and hence there are $x,y,z\in\N$
 with $\gcd(x,y,z)=1$ such that $4n+1=(dx)^2+(dy)^2+(dz)^2$ with $dx-dy,dx-dz,dy-dz\not=\pm1$.
 When $4n+1$ is squarefree, by Lemma 2.3 we have $h(-4n-1)\ls 4$ if and only if $4n+1$ belongs to the set
 $$S:=\{13,\ 17,\ 21,\ 33,\ 37,\ 57,\ 73,\ 85,\ 93,\ 97,\ 133,\ 177,\ 193,\ 253\}.$$
 If $4n+1\in S$ then we can easily write $4n+1$ as $x^2+y^2+z^2$ with $x,y,z\in\N$, $2\nmid z$ and $y-z\not=\pm1$.
 (For example, $4\times63+1=253=10^2+12^2+3^2$.)
 When $4n+1$ is squarefree with $4n+1\not\in S$, by Lemma 2.1 we have
 $R_3(4n+1)\gs 12h(-4n-1)>12\times4=48$
 and hence $4n+1=x^2+y^2+z^2$ for some $x,y,z\in\N$ with $2\nmid z$ and $y-z\not=\pm1$.

 In view of the above, we can write $4n+1$ as $(2x)^2+(2y)^2+(2z+1)^2$ with $x,y,z\in\N$ and $2y-(2z+1)\not=\pm1$.
 It follows that $n=x^2+y^2+2T_z=x^2+T_{y+z}+T_{z-y}$ with $z-y\not\in\{0,-1\}$. Note that $T_{z-y}\not=0$.
 If $y+z\ls 5$, then $T_{y+z}+T_{z-y}$ belongs to the set
 $$T:=\{T_i+T_j:\, i,j=0,\ldots,5\}
 =\{0,1,3,4,6,7,9,\ldots,13,15,16,18,20,21,25,30\}.$$
If $n=x^2+t>66=6^2+30$ with $x\in\N$ and $t\in T$, then  $x>6$ and hence by (c) we can write $x^2+t$ as $a^2+T_b+T_c$ with $a,b,c\in\N$, $b\gs c>0$ and $b\gs6$.

 Let $k\in\{1,2,3,4,5,6\}$. It is easy to verify that each $n=T_k+1,T_k+2,\ldots,66$ can be written as $x^2+T_y+T_z$ with $x,y,z\in\N$, $y\gs k$ and $z\gs1$.
 Thus, for any $n\in\N$ we can write $n+T_k+1$ as $x^2+T_{y+k}+T_{z+1}$ with $x,y,z\in\N$ and hence
 $$n=x^2+T_{y+k}-T_k+T_{z+1}-T_1=x^2+\f{y(y+2k+1)}2+\f{z(z+3)}2.$$
 This proves the universality of $(2,0,1,2k+1,1,3)$ over $\N$.

(e) As conjectured by Sun [S07] and proved in [OS],
 any positive integer can be written as the sum of a square, an odd square and a triangular number.
 So, for any $m\in\Z^+$ there are $a,b,c\in\N$ with $a$ odd such that $T_m=a^2+b^2+T_c$ and hence
 $2T_m=(a+b)^2+|a-b|^2+2T_c$ with $a+b\gs a>0$.
 If $n\in\Z^+$ is not twice a triangular number, then $4n+1$ is not a square, hence by the Gauss-Legendre theorem
 there are  $x,y,z\in\N$ with $\max\{x,y\}>0$ such that
 $4n+1=(2x)^2+(2y)^2+(2z+1)^2$ and thus
 $n=x^2+y^2+2T_z$.

 By the above, for any $n\in\N$ there are $x,y,z\in\N$
 such that  $n+1=(x+1)^2+y^2+2T_z=x(x+2)+1+T_{y+z}+T_{z-y}$.
 So $(2,4,2,2,2,0)$ and $(2,4,1,1,1,1)$ are universal over $\N$.

 Clearly, $4x^2+T_y+z(z+3)/2$ is universal over $\N$
if and only if any positive integer can be written as $4x^2+T_y+T_z$ with $x,y,z\in\N$ and $\max\{y,z\}>0$.
By Sun [S07, Theorem 1(i)], we can write any $n\in\Z^+$ in the form $(2x)^2+T_y+T_z$ with $x,y,z\in\N$.
If $y=z=0$, then $n=4x^2=4\times0^2+T_{2x-1}+T_{2x}$. Thus $4x^2+T_y+z(z+3)/2$ (or the tuple $(8,0,1,3,1,1)$) is universal over $\N$.

 (f) Let $n>3$ be an integer. If $4n+1=m^2$ with $m\in\{4,5,6,\ldots\}$, then by Lemma 2.5 we can write $4n+1=x^2+y^2+z^2$ with $x,y,z\in\N$, $2\nmid z$ and $\max\{x,y\}>0$,
 hence $\{x,y\}\not\se\{0,\pm2\}$ since $m^2-z^2\not\in\{2^2,2^2+2^2\}$.
 If $4n+1$ is not a square but $d^2\mid 4n+1$ for some integer $d>1$, then by Lemma 2.1 there are $x,y,z\in\N$ with $2\nmid z$ such that
 $4n+1=(dx)^2+(dy)^2+(dz)^2$ with $\max\{dx,dy\}\gs d>2$.

 Now suppose that $4n+1$ is squarefree. If $n=9$, then $4n+1=0^2+6^2+1^2$. If $n\not=9$, then
 by Lemmas 2.1 and 2.3 we have $R_3(4n+1)=12h(-4n-1)>12\times2=24$
 and hence the equation $x^2+y^2+z^2=4n+1$ has integral solutions with $2\nmid z$ and $\{x,y\}\not\se\{0,\pm2\}$.
 (As $4n+1>16$, there is at most one value of $\da\in\{0,1,2\}$ with $4n+1-\da2^2$ a square.)

 By the above, for any integer $n\gs2$ there are $x,y,z\in\N$ such that
 $4(n+2)+1=(2(x+2))^2+(2y)^2+(2z+1)^2$
 and hence $n=x(x+4)+y^2+2T_z=x(x+4)+T_{y+z}+T_{z-y}$.
Note also that $n=0(0+4)+n^2+2T_0=0(0+4)+T_n+T_0$ for each $n=0,1$.
Therefore both $(2,8,2,2,2,0)$ and $(2,8,1,1,1,1)$ are universal tuples over $\N$.

Let $n\in\N$. As mentioned in the last paragraph, there are $x,y,z\in\N$ such that
$$4(2n+5)+1=(2x+4)^2+(2y)^2+(2z+1)^2=2(x+y+2)^2+2(x-y+2)^2+(2z+1)^2.$$
Since $(x+y+2)^2+(x-y+2)^2\eq 2\pmod 4$, there are $u,v\in\N$ such that $x+y+2=2u+3$ and $|x-y+2|=2v+1$.
Thus $8n+21=2(2u+3)^2+2(2v+1)^2+(2z+1)^2$ and hence
$n=u(u+3)+v(v+1)+T_z$. This proves the universality of $(2,6,2,2,1,1)$ over $\N$.

(g) Let $n\gs2$ be an integer. By Lemma 2.5, there are $x,y,z\in\N$ with $x\gs y\gs z$ and $y>0$ such that $n=T_x+T_y+T_z$.
If $\{x-z,y-z\}\se\{1,3,5\}$, then $x-y\not\in\{1,3,5\}$ since $x\eq y\pmod2$, and $n\ls T_5+T_5+T_0=30$ if $z=0$.

Let $n>31$ be an integer. By the last paragraph, $n=T_x+T_y+T_z$ for some $x,y\in\Z^+$ and $z\in\{0,\ldots,y\}$ with $y-z\not\in\{1,3,5\}$.
If $y+z+1\in\{1,3,5\}$, then $(y,z)\in\{(1,1), (2,0), (2,2),(3,1),(4,0)\}$
and hence $T_y+T_z\ls T_4=10$, thus $x>6$ (as $n>T_6+10=31$), and $x-z,x+z+1\not\in\{1,3,5\}$ unless $(y,z)=(2,2)$ and $x=7$
in which case $n=T_7+T_2+T_2=T_7+T_3+T_0$ with $7-0,7+0+1\not\in\{1,3,5\}$.
Therefore, we can always write $n$ as $T_u+T_v+T_w$ with $u,v\in\Z^+$, $w\in\{0,\ldots,v\}$ and $v-w,v+w+1\not\in\{1,3,5\}$.
It follows that
$$T_v+T_w=\f{(v+w+1)^2+(v-w)^2-1}4=\f{(2r)^2+(2s+1)^2-1}4=r^2+s(s+1)$$
for some $r,s\in\N$ with $s>2$. Hence $n=T_u+r^2+s(s+1)$.
For each $k\in\{0,1,2,3\}$, clearly $t=s-k\in\N$ and
$$n-k(k+1)=T_u+r^2+(t+k)(t+k+1)-k(k+1)=T_u+r^2+t(t+2k+1).$$
Note also that $T_u=T_{a+1}=a(a+3)/2+1$ with $a=u-1\in\N$.

Let $k\in\{1,2,3\}$. By the above, any integer $n>30-k(k-1)$ can be written $a(a+3)/2+r^2+t(t+2(k-1)+1)$ with $a,r,t\in\N$.
We can easily see that each $n=0,1,\ldots,30-k(k-1)$ also can be written as $x(x+3)/2+y^2+z(z+2k-1)$ with $x,y,z\in\N$.
Thus the tuple $(2,4k-2,2,0,1,3)$ is universal over $\N$.
Similarly,  $(2,4k+2,2,0,1,1)$ is also a universal tuple over $\N$.

(h) Let $k\in\{1,2,3\}$. It is easy to verify that every $n=k^2,\ldots,33$ can be written as $T_x+y^2+z^2$ with $x,y,z\in\N$ and $\max\{y,z\}\gs k$.

Now let $n\in\N$ with $n>33$. By [OS], $n=T_x+y^2+z^2$ for some $x,y,z\in\N$ with $y$ odd.
Assume that $\max\{y,z\}<k$. Then $y=1$ and $z\ls k-1$, hence
$$r:=n-T_x\in\{s^2+1:\ s=0,\ldots,k-1\}\se\{s^2+1:\ s=0,1,2\}=\{1,2,5\}.$$
As $n>33$, we have $T_x>28$ and hence $x>7$. Since
$T_x<n=T_x+r\ls T_x+5<T_{x+1}$,
$n$ is not a triangular number. By [S07, Theorem 1(ii)], there are $a,b,c,u,v,w\in\N$ with $a\not\eq b\pmod 2$
and $u\eq v\pmod 2$ such that
$n=a^2+b^2+T_c=u^2+v^2+T_w.$
Suppose that $\max\{a,b\}<k\ls 3$ and also $\max\{u,v\}<k\ls 3$.
Then $a^2+b^2\in\{1^2+0^2,1^2+2^2\}$ and $u^2+v^2\in\{0^2+0^2,0^2+2^2,1^2+1^2,2^2+2^2\}$.
Note that $T_{x+1}>n\gs T_w=n-(u^2+v^2)\gs T_x+1-8>T_{x-1}.$
So we have $w=x$ and hence $r=u^2+v^2\in\{0,2,4,8\}$. Similarly, $c=x$ and hence $r=a^2+b^2\in\{1,5\}$.
Thus we get a contradiction since $\{0,2,4,8\}\cap\{1,5\}=\emptyset$.

In view of the above, $T_x+y^2+z(z+2k)=T_x+y^2+(z+k)^2-k^2$ is universal over $\N$, i.e., the tuple $(2,4k,2,0,1,1)$ is universal over $\N$.

(i) Let $n\in\N$. By Lemma 2.1 we have $R_3(8(n+1)+3)>0$ and hence $n+1=T_x+T_y+T_z$ for some $x,y,z\in\N$ with $y>z$. Clearly,
$4(T_y+T_z)+1=(y+z+1)^2+(y-z)^2=(2(u+1))^2+(2v+1)^2$ for some $u,v\in\N$, and hence $n=T_x+u(u+2)+v(v+1)$. Therefore
$(2,4,2,2,1,1)$ is universal over $\N$.

(j) Let $k\in\{1,2,3\}$. Clearly, $(2,4k,2,0,1,3)$ is universal over $\N$ if and only if any integer $n\gs k^2+1$ can be written as $x^2+y^2+T_z$ with $x,y,z\in\N$, $y\gs k$ and $z\gs1$.
By Lemma 2.8 and the equalities $T_2=1^2+1^2+T_1$ and $T_3=1^2+2^2+T_1$, any triangular number greater than $k^2$ can be written as $x^2+y^2+T_z$ with $x,y,z\in\N$, $y\gs k$ and $z\gs1$.

Now we fix an integer $n\gs k^2+1$ which is not a triangular number. Then $8n+1$ is not a square. By [S07, Theorem 1(ii)],
$$r_0(n)=r_1(n)=\f14|\{(x,y,z)\in\Z^3:\ x^2+y^2+T_z=n\}|,$$
where
$$ r_{\da}(n):=|\{(x,y,z)\in\Z\times\Z\times\N:\ x^2+y^2+T_z=n\ \t{and}\ x-y\eq \da\ (\t{mod}\ 2)\}|$$
for $\da=0,1$.
For $x,y,z\in\Z$, clearly
$$n=x^2+y^2+T_z\iff 8n+1=(2x+2y)^2+(2x-2y)^2+(2z+1)^2.$$
If $8n+1=u^2+v^2+w^2$ with $u,v,w\in\N$ and $2\nmid w$, then $8\mid u^2+v^2$, hence $2\mid u$, $2\mid v$ and $u/2\eq v/2\pmod2$.
Thus, with the help of the Gauss-Legendre theorem we have
$$12r_0(n)=12r_1(n)=|\{(x,y,z)\in\Z^3:\ x^2+y^2+z^2=8n+1\}|>0.\tag3.1$$

If $n-0^2-1^2$ and $n-2^2-1^2$ are both triangular numbers, then $n=11$.
If at least two of $n-1^2-1^2,n-0^2-2^2,n-2^2-2^2$ are triangular numbers, then $n\in\{14,23\}$.
When $n\in\{11,14,23\}$, we can easily verify that $n=x^2+y^2+T_z$ for some $x,y,z\in\N$ with $y\gs k$ and $z\gs1$.
Let $n\not=11,14,23$. If $r_0(n)=r_1(n)>4(k-1)$, or $r_0(n)=r_1(n)>4$ and $n-5\not\in\{T_m:\ m\in\N\}$, then $n=a^2+b^2+T_c=u^2+v^2+T_w$
for some $a,b,c,u,v,w\in\N$ with $a\not\eq b\pmod 2$, $u\eq v\pmod2$, $\max\{a,b\}\gs k$ and $\max\{u,v\}\gs k$,
hence $T_c\not\eq T_w\pmod2$ and thus $n=x^2+y^2+T_z$ for some $x,y,z\in\N$ with $y\gs k$ and $z\gs1$.

The above arguments with $k=1$ yield that any integer $n\gs 1^2+1$ can be written as $x^2+(y+1)^2+T_{z+1}$ with $x,y,z\in\N$.
Thus $(2,0,2,0,1,3)$ and $(2,4,2,0,1,3)$ are universal tuples over $\N$.

Now we assume that $k\in\{2,3\}$. Write $8n+1=d^2q$ with $d,q\in\Z^+$ and $q$ squarefree. Then $q>1$ and $q\eq1\pmod8$. By Lemma 2.3, $h(-8n-1)=h(-q)\gs4$, and the equality holds only when $q$ belongs to the set
$Q=\{17,\ 33,\ 57,\ 73,\ 97,\ 177,\ 193\}.$
By Lemma 2.1, if $d>1$ then
$$|\{(x,y,z)\in\Z^3:\ x^2+y^2+z^2=8n+1\}|\gs R_3(8n+1)+R_3(q)\gs 24h(-q)$$
and hence $r_0(n)=r_1(n)\gs 2h(-q)\gs8$.
If $n=1^2+2^2+T_m$ for some $m\in\N$, then $(2m+1)^2+40=8n+1=d^2q$ and hence $q\not\in Q$
since $(\f{-40}p)=-1$ for all $p\in Q$ with $(\f{\cdot}p)$ the Jacobi symbol,
therefore $r_0(n)=r_1(n)\gs 2h(-q)>8$ when $d>1$.

Now we handle the case $d=1$. If $8n+1=q\in Q$, then $n\in\{2,4,7,9,12,22,24\}$. Recall that $n\gs k^2+1$.
Clearly,
$$7=0^2+2^2+T_2,\ 9=2^2+2^2+T_1,\ 12=0^2+3^2+T_2,\ 22=0^2+4^2+T_3,\ 24=0^2+3^2+T_5.$$
Now assume that $8n+1=q\not\in Q$. Then $h(-q)>4$. If $h(-q)\ls 8$, then by Lemma 2.3(v), $n$ belongs to the set
$$\{5,\, 8,\,13,\,14,\,17,\,18,\,27,\,33,\,34,\,39,\,42,\,43,\,48,\,57,\,66,\,69,\,87,\,99\}$$
and we can directly verify that $n=x^2+y^2+T_z$ for some $x,y,z\in\N$ with $y\gs k$ and $z\gs1$.
When $h(-q)>8$, by (3.1) and Lemma 2.1 we have
$r_0(n)=r_1(n)=R_3(8n+1)/12\gs h(-8n-1)=h(-q)>8.$
Thus $(2,4k,2,0,1,3)$ is indeed universal over $\N$.

(k) Let $k\in\{1,2,3\}$. We want to prove that $2x^2+T_y+z(z+2k+1)/2=2x^2+T_y+T_{z+k}-T_k$ is universal over $\N$ (i.e., $(4,0,1,2k+1,1,1)$ is a universal tuple over $\N$).
It is easy to check that each $n=T_k,T_k+1,\ldots,77$ can be written as $x^2+T_y+T_z$ with $x,y,z\in\N$ and $\max\{y,z\}\gs k$.

Now let $n\in\Z^+$ with $n\gs 78$. If $8n+2=x^2+a=y^2+b$ with $x,y\in\N$, $a<b$ and
$$a,b\in A:=\{w^2+z^2:\ w,z\in\{1,3,5\}\}=\{2,10,18,26,34,50\},$$
then $x>y\gs \sqrt{8\times 78+2-50}=24$ and hence
$x^2+a>y^2+2y+a\gs y^2+48+a\gs y^2+b$.
So $8n+2$ can be written as $x^2+a$ with $a\in A$ in at most one way. Therefore,
the equation $8n+2=x^2+y^2+z^2$ has at most $3!\times 2^3=48$ integral solutions with the two odd numbers among $x,y,z$ in the set $\{\pm1,\pm3,\pm5\}$.

Write $8n+2=d^2q$ with $d,q\in\Z^+$ and $q$ squarefree. Then $8\mid(q-2)$. By Lemma 2.3, $h(-q)\ls 4$ only when $q$ is among
$2,\, 10,\ 34,\ 42,\ 82,\ 130.$
Clearly, $82d^2=0^2+d^2+(9d)^2$ and $130d^2=0^2+(7d)^2+(9d)^2$.
If $h(-q)>4$, then by Lemma 2.1 we have
$$|\{(x,y,z)\in\Z^3:\ x^2+y^2+z^2=8n+2\}|\gs R_3(q)=12h(-q)>48$$
and hence $8n+2=x^2+y^2+z^2$ for some $x,y,z\in\N$ with the two odd numbers of $x,y,z$ not all in $\{1,3,5\}$.
As $R_3((8n+2)/d^2)>0$ by Lemma 2.1, there are $x,y,z\in\N$ with $\gcd(x,y,z)=1$ such that $8n+2=(dx)^2+(dy)^2+(dz)^2$
with two of the odd numbers among $dx,dy,dz$ at least $d$.
If $d\in\{1,3\}$, then $q=(8n+2)/d^2\gs (8\times78+2)/9=626/9>69$. If $d=5$ then $q\gs 626/25>25$.
Note that $5^2\times34=0^2+15^2+25^2$ and $5^2\times42=20^2+5^2+25^2$.
Therefore, we always can write $8n+2=(2w)^2+(2y+1)^2+(2z+1)^2$ with $w,y,z\in\N$ and $\max\{y,z\}\gs3\gs k$.
Clearly, $w=2x$ for some $x\in\N$, and hence $n=2x^2+T_y+T_z$. We are done.

(l) To prove that $(3,9,2,0,1,3)$, $(3,3,2,0,1,3)$ and $(3,9,2,0,1,1)$ are universal tuples over $\N$,
we only need to show that any integer $n\gs 4$ can be written as $x^2+T_y+3T_z$ with $x\in\N$ and $y,z\in\Z^+$
which can be easily verified for all $n=4,5,\ldots,45$.

Now, we fix an integer $n\gs 46$. By the Gauss-Legendre theorem, $12n+6=x^2+y^2+z^2$ for some $x,y,z\in\Z$ with $2\mid x$ and $2\nmid yz$. Note that $12n+6>12\times46+6>9^2+(9+3)^2+(9+9)^2$.
It is easy to see that the equation $12n+6=x^2+y^2+z^2$ has at most $2!\times 2^3=16$ integral solutions with $x$ even, and $y,z\in\{\pm(x-3),\pm(x+3)\}$ or
$$\{y,z\}\in\{\{\ve_1(x+3\ve_0),\ve_2(x+9\ve_0)\}:\ \ve_0,\ve_1,\ve_2\in\{\pm1\}\}.\tag3.2$$

Write $12n+6=d^2q$ with $d,q\in\Z^+$ and $q$ squarefree. Obviously $q\eq 2\pmod4$. If $d>3$, then by Lemma 2.1 we have $12n+6=x^2+y^2+z^2$
for some $x,y,z\in\N$ with $2\mid x$ and $\gcd(x,y,z)=d>3$, hence $\{y,z\}\not\se\{\pm(x-3),\pm(x+3)\}$ and (3.2) fails.
If $h(-q)>4$, then by Lemma 2.3 we have
$R_3(12n+6)=12h(-12n-6)=12h(-q)>48$
and hence the equation $12n+6=x^2+y^2+z^2$ has more than 16 solutions with $x$ even.
If $h(q)\ls4$, then by Lemma 2.3 the number $q$ belongs to the set
$$E=\{2,\,6,10,\,14,\,\,22,\,30,\,34,\,42,\,46,\,58,\,70,\,78,\,82,\,102,\,130,\,142,\,190\}.$$
If $d=1$, then $q=12n+6>12\times 46+6>190$ and hence $h(-q)>4$.
If $12n+6=3^2q$ with $q\in E$, then we can verify that the equation $12n+6=x^2+y^2+z^2$ has solutions with $x$ even, $\{y,z\}\not\se\{\pm(x-3),\pm(x+3)\}$ and $(3.2)$ invalid.

By the above, there are $x,y,z\in\Z$ with $2\mid x$, $2\nmid yz$ and $\{y,z\}\not\se\{\pm(x-3),\pm(x+3)\}$ such that (3.2) fails.
As $x^2+y^2+z^2\eq0\pmod3$, either $x\eq y\eq z\eq0\pmod3$ or $3\nmid xyz$.
Without loss of generality, we may assume that $x\eq y\eq z\pmod 3$.
Recall Jacobi's identity
$$3(x^2+y^2+z^2)=(x+y+z)^2+2\l(\f{x+y}2-z\r)^2+6\l(\f{x-y}2\r)^2\tag3.3$$
which can be verified directly. Clearly, $x+y+z=6w$ for some $w\in\Z$,
and $x+y-2z=6u+3$ and $x-y=6v+3$ for some $u,v\in\Z$. Thus, by (3.3) we have
$$36n+18=3(x^2+y^2+z^2)=(6w)^2+2\l(\f{6u+3}2\r)^2+6\l(\f{6v+3}2\r)^2$$
and hence $n=w^2+T_u+3T_v$. If $x-y\not=\pm3$, then $v\not=0,-1$ and hence $T_v\not=0$.
If $x+y-2z\not=\pm3$, then $u\not=0,-1$ and hence $T_u\not=0$.
So we are done if $\{x-y,x+y-2z\}\cap\{\pm3\}=\em$. Due to the symmetry of $y$ and $z$, we are also done if $\{x-z,x+z-2y\}\cap\{\pm3\}=\em$.

If $\{x+y-2z,x+z-2y\}\se\{\pm3\}$, then $3(y-z)=x+y-2z-(x+z-2y)\in\{0,\pm6\}$, hence $y=z$ (since $y\eq z\pmod3$) and $x-y=x-z\in\{\pm3\}$,
which contradicts that $\{y,z\}\not\se\{\pm(x-3),\pm(x+3)\}$. If $x-y,x-z\not=\pm3$, then we are done since $x+y-2z\not=\pm3$ or $x+z-2y\not=\pm3$.

Now we consider the remaining case in which exactly one of $|x-y|$ and $|x-z|$ is 3. Without loss of generality, we assume that
$x-y\not=\pm3$ and $x-z\in\{\pm3\}$. Thus $y\not=z$. We are done if $x+y-2z\not=\pm3$. Suppose that $x+y-2z=(x-z)+(y-z)\in\{\pm3\}$.
Then $y-z=-2(x-z)\in\{\pm6\}$ and $x-y=x-z-(y-z)=3(x-z)\in\{\pm9\}$. So $(y,z)=(x+9,x+3)$ or $(x-9,x-3)$, which contradicts that (3.2) fails.
\smallskip

In view of the above, we have completed the proof of Theorem 1.1. \qed

\heading{4. Proofs of Theorems 1.4 and 1.7}\endheading

\proclaim{Lemma 4.1} Let $a,b,c,d\in\Z$ with $a\gs c\gs1$,  $b>-a$, $d>-c$, $a\eq b\pmod 2$ and $c\eq d\pmod2$.
Then $\{1,\ldots,18\}\not\se\{\psi_{a,b}(x)+\psi_{c,d}(y):\ x,y\in\N\}$.
\endproclaim
\Proof. It is easy to see that neither $\{x(ax+b)/2:\, x\in\N\}$ nor $\{y(cy+d)/2:\,y\in\N\}$ contains $\{1,2\}$.
So $\{1,2\}\not\se\{\psi_{a,b}(x)+\psi_{c,d}(y):\ x,y\in\N\}$
if $\psi_{a,b}(1)=(a+b)/2>2$ or $\psi_{c,d}(1)=(c+d)/2>2$.

Below we suppose that $a+b\ls 4$ and $c+d\ls 4$.
In the case $ac<212$, via a computer we find that one of $1,\ldots,9$ cannot be written as $\psi_{a,b}(x)+\psi_{c,d}(y)$ with $x,y\in\N$.

Now we assume that $ac\gs 212$. Then $1/a+1/c\ls 213/212$.

Fix a positive integer $N$. For $x\in\Z$, it is easy to see that
$$\f {x(ax+b)}2\ls N\iff -\f{\sqrt{8aN+b^2}+b}{2a}\ls x\ls \f{\sqrt{8aN+b^2}-b}{2a}.\tag4.1$$
As $-a<b\ls a$ or $b>0$, we have
$$\l|\l\{x\in\N:\ \psi_{a,b}(x)\ls N\r\}\r|\ls1+\f{\sqrt{8aN+b^2}-b}{2a}<\f32+\sqrt{\f{2N}a+\f14};$$
Similarly, $|\{x\in\N:\ \psi_{c,d}(x)\ls N\}|<3/2+\sqrt{2N/c+1/4}.$
Note that
$$\sqrt u+\sqrt v\ls\sqrt{2u+2v}\quad\t{for all}\ u,v\gs0.\tag4.2$$
Therefore
$$\align&\l|\l\{\psi_{a,b}(x)+\psi_{c,d}(y):\ x,y\in\N\r\}\cap[0,N]\r|
\\<&\l(\f32+\sqrt{\f{2N}a+\f14}\r)\l(\f32+\sqrt{\f{2N}c+\f14}\r)
\\\ls&\f94+\sqrt{\f{4N^2}{ac}+\f N2\l(\f1a+\f1c\r)+\f1{16}}+\f32\sqrt{4N\l(\f1a+\f1c\r)+1}
\\\ls&f(N):=\f94+\sqrt{\f{4N^2}{212}+\f N2\cdot\f{213}{212}+\f1{16}}+\f32\sqrt{4N\times\f{213}{212}+1}.
\endalign$$

Now, take $N=18$. Then
$$\l|\l\{\psi_{a,b}(x)+\psi_{c,d}(y):\ x,y\in\N\r\}\cap[0,N]\r|\ls f(N)<1+N.$$
So one of $1,\ldots,N$ cannot be written as $\psi_{a,b}(x)+\psi_{c,d}(y)$ with $x,y\in\N$.
\qed

\medskip
\noindent{\it Proof of Theorem 1.4}. In view of Lemma 4.1, for certain $m\in\{1,\ldots,18\}$ we can write $m=\psi_{a,b}(x)+\psi_{c,d}(y)+\psi_{e,f}(z)$ with $x,y,z\in\N$ and $x>0$, and hence
$(a+b)/2\ls\psi_{a,b}(x)\ls m\ls 18.$
Similarly, $(c+d)/2\ls18$ and $(e+f)/2\ls 18$. Therefore
$$b\ls 36-a,\ d\ls 36-c\ \t{and}\ f\ls 36-e.\tag4.3$$

If $ce\gs2000$ and $N=64$, then $1/c+1/e\ls 2001/2000$ and hence
$$\align
&\l|\l\{\psi_{c,d}(y)+\psi_{e,f}(z):\ y,z\in\Z\r\}\cap[0,N]\r|
\\<&\f 94+\sqrt{\f{4N^2}{2000}+\f N2\cdot\f{2001}{2000}+\f1{16}}+\f32\sqrt{4N\times\f{2001}{2000}+1}\ls 33
\endalign$$
by the proof of Lemma 4.1, thus
$$\l|\l\{\psi_{a,b}(x)+\psi_{c,d}(y)+\psi_{e,f}(z):\ x\in\{0,1\}\ \&\ y,z\in\Z\r\}\cap[0,N]\r|$$
is at most $2\times32=N$. So, in the case $ce\gs2000$, we can write certain $n\in\{0,\ldots,64\}$
as $\psi_{a,b}(x)+\psi_{c,d}(y)+\psi_{e,f}(z)$
with $x,y,z\in\N$ and $x\gs2$, hence
$a+2\ls2a+b\ls \psi_{a,b}(x)\ls n\ls 64$
and thus $e\ls c\ls a\ls 62$.

Now we consider the case $ce<2000$.  In view of (4.3), via a computer we find that for each $i=1,\ldots,18$ there is an integer $n_i\in[i,58]$ such that
$$\{n_i,n_i-i\}\cap\l\{\psi_{c,d}(y)+\psi_{e,f}(z):\ y,z\in\N\r\}=\em.$$
(We note that $\{58,58-17\}\cap\{y^2+p_3(z):\ y,z\in\N\}=\em$.)
For $i=(a+b)/2\ls 18$, we can write $n_i$ as $\psi_{a,b}(x)+\psi_{c,d}(y)+\psi_{e,f}(z)$ with $x,y,z\in\N$ and $x>1$, thus
$a+2\ls 2a+b\ls\psi_{a,b}(x)\ls n_i\ls 58$
and hence $e\ls c\ls a\ls 56$.

By the above, either $ce\gs2000$ and $e\ls c\ls a\ls 62$, or $ce<2000$ and $e\ls c\ls a\ls58$.
In view of this and (4.3), via a computer we find that if every $n=0,\ldots,10^5$ can be written as $\psi_{a,b}(x)+\psi_{c,d}(y)+\psi_{e,f}(z)$ with $x,y,z\in\Z$
then the tuple $(a,b,c,d,e,f)$
must be among the 56+10 tuples listed in Theorem 1.1 and Conjecture 1.2 if $a\mid b$, $c\mid d$ and $e\mid f$, or among the 407 tuples listed the Appendix if $a\nmid b$ or $c\nmid d$
or $e\nmid f$. This concludes the proof of Theorem 1.4. \qed

\proclaim{Lemma 4.2} Let $a,b,c,d\in\N$ with $a\gs c\gs2$, $a>b$, $c>d$, $a\eq b\pmod 2$ and $c\eq d\pmod2$.
Then one of $1,\ldots,28$ cannot be written as $\psi_{a,b}(x)+\psi_{c,d}(y)$ with $x,y\in\Z$.
\endproclaim
\Proof. If $ac<1000$, then $\max\{a,b,c,d\}=a<500$, and hence we may use a computer to get that
one of $1,\ldots,21$ cannot be written as $\psi_{a,b}(x)+\psi_{c,d}(y)$ with $x,y\in\Z$. Note that $21\not=x(7x+1)/2+y(3y+1)/2$ for all $x,y\in\Z$.

Below we assume that $ac\gs 1000$. It is easy to see that $1/a+1/c\ls251/500$.

Let $N$ be any positive integer. In view of (4.1),
$$\l|\l\{x\in\Z:\ \psi_{a,b}(x)\ls N\r\}\r|-1\ls \f{\sqrt{8aN+b^2}-b}{2a}-\l(-\f{\sqrt{8aN+b^2}+b}{2a}\r)<\sqrt{\f{8N}a+1}.$$
Similarly,
$\l|\l\{y\in\Z:\ \psi_{c,d}(y)\ls N\r\}\r|<1+\sqrt{8N/c+1}.$
With the help of (4.2), we have
$$\align&\l|\l\{\psi_{a,b}(x)+\psi_{c,d}(y):\ x,y\in\Z\r\}\cap[0,N]\r|
\\<&\l(1+\sqrt{\f{8N}a+1}\r)\l(1+\sqrt{\f{8N}c+1}\r)
\\\ls&1+\sqrt{\f{64N^2}{ac}+8N\l(\f1a+\f1c\r)+1}+\sqrt{16N\l(\f1a+\f1c\r)+4}
\\\ls&g(N):=1+\sqrt{\f{64N^2}{1000}+8N\f{251}{500}+1}+\sqrt{16N\f{251}{500}+4}.
\endalign$$

Now, take $N=28$. Then
$$\l|\l\{\psi_{a,b}(x)+\psi_{c,d}(y):\ x,y\in\Z\r\}\cap[0,N]\r|<g(N)\ls 1+N.$$
Therefore, one of $1,\ldots,N$ cannot be written as $\psi_{a,b}(x)+\psi_{c,d}(y)$ with $x,y\in\Z$.
\qed

\medskip
\noindent{\it Proof of Theorem 1.7}. For any integer $x$ with $|x|\gs2$, we have
$$\psi_{a,b}(x)\gs\f{|x|(a|x|-b)}2\gs a|x|-b\gs 2a-b>\psi_{a,b}(1)=\f{a+b}2.$$

In view of Lemma 4.2, for certain $m\in\{1,\ldots,28\}$ we can write $m=\psi_{a,b}(x)+\psi_{c,d}(y)+\psi_{e,f}(z)$ with $x,y,z\in\Z$ and $x\not=0$, and hence
$(a-b)/2=\psi_{a,b}(-1)\ls \psi_{a,b}(x)\ls m\ls 28.$
Similarly, $(c-d)/2\ls28$ and $(e-f)/2\ls 28$.

If $ce\gs1000$ and $N=190$, then
$$\l|\l\{\psi_{c,d}(y)+\psi_{e,f}(z):\ y,z\in\Z\r\}\cap[0,N]\r|<g(N)\ls 96$$
by the proof of Lemma 4.2, and hence
$$\l|\l\{\psi_{a,b}(x)+\psi_{c,d}(y)+\psi_{e,f}(z):\ x\in\{0,-1\}\ \&\ y,z\in\Z\r\}\cap[0,N]\r|$$
is at most $2\times95=N$. So, in the case $ce\gs1000$, we can write certain $n\in\{0,1,\ldots,190\}$
as $\psi_{a,b}(x)+\psi_{c,d}(y)+\psi_{e,f}(z)$
with $x,y,z\in\Z$ and $x\not=0,-1$, hence
$(a+b)/2\ls \psi_{a,b}(x)\ls n\ls 190$
and thus $a=(a-b)/2+(a+b)/2\ls 28+190=218.$

Now we consider the case $ce<1000$.  Via a computer we find that for each $i=1,\ldots,28$ there is an integer $n_i\in[i,157]$ such that
$$\{n_i,n_i-i\}\cap\l\{\psi_{c,d}(y)+\psi_{e,f}(z):\ y,z\in\Z\r\}=\em.$$
For $i=(a-b)/2\ls28$, we can write $n_i$ as $\psi_{a,b}(x)+\psi_{c,d}(y)+\psi_{e,f}(z)$ with $x,y,z\in\Z$ and $x\not=0,-1$, thus
$(a+b)/2\ls\psi_{a,b}(x)\ls n_i\ls 157$
and hence $a=(a-b)/2+(a+b)/2\ls 28+157=185.$

By the above, either $ce\gs1000$ and $e\ls c\ls a\ls 218$, or $ce<1000$ and $e\ls c\ls a\ls 185$.
Via a computer we find that if each $n=0,\ldots,10^5$ can be written as $\psi_{a,b}(x)+\psi_{c,d}(y)+\psi_{e,f}(z)$ with $x,y,z\in\Z$
then the tuple $(a,b,c,d,e,f)$ must be among the 12082 tuples listed in [S17a].
This completes the proof of Theorem 1.7. \qed

\heading{5. Proofs of Theorems 1.9 and 1.11}\endheading

\medskip
\noindent{\it Proof of Theorem 1.9}. Let $n$ be any nonnegative integer.

(i) $n=x(ax+b)+y(ay+a-b)$ if and only if $4an+b^2+(a-b)^2$ coincides with
$$(a(2x)+b)^2+(a(2y)+a-b)^2=(a(-2y-1)+b)^2+(a(-2x-1)+a-b)^2.$$
Therefore
$$\align &n\in\{x(ax+b)+y(ay+a-b):\ x,y\in\Z\}
\\\iff& 4an+(a-b)^2+b^2\in\{(au+b)^2+(av+a-b)^2:\ u,v\in\Z\ \&\ 2\mid u-v\}
\\\iff&4an+(a-b)^2+b^2\in\{(a(x-y)+b)^2+(a(x+y)+a-b)^2:\ x,y\in\Z\}
\\\iff&n\in\{aT_x+\psi_{a,a-2b}(y):\ x,y\in\Z\}.
\endalign$$
This proves (1.6).

(ii) Obverse that $2(x^2+4T_y)+1=2x^2+(2y+1)^2$ and
$$\align &2n+1\in\{u^2+2v^2:\ u,v\in\Z\}
\\\iff&2n+1\in\l\{\l(\f{x-y}3+y\r)^2+2\l(\f{x-y}3\r)^2:\ x,y\in\Z\ \&\ 3\mid x-y\r\}
\\\iff&6n+3\in\{x^2+2y^2:\ x,y\in\Z\ \&\ 3\mid x-y\}
\\\iff&6n+3\in\{x^2+2y^2:\ x,y\in\Z\}.
\\\iff&6n+3\in\{x^2+2y^2:\ x,y\in\Z,\ 2\nmid x\ \&\ 3\nmid xy\}
\\&\qquad\qquad\qquad\ (\t{by [JP, p.\,173] or [S15, Lemma 2.1]})
\\\iff&6n+3\in\{(6x-1)^2+2(3y-1)^2:\ x,y\in\Z\}
\\\iff& n=4p_5(x)+p_8(y)\ \t{for some}\ x,y\in\Z.
\endalign$$
So (1.9) holds. Similarly,
$$\align n\in\{x^2+T_y:\ x,y\in\Z\}
\iff&8n+1\in\{u^2+2v^2:\ u,v\in\Z\}
\\\iff&24n+3\in\{u^2+2v^2:\ u,v\in\Z,\ 2\nmid uv\ \&\ 3\nmid uv\}
\\\iff&24n+3\in\{(6x-1)^2+2(6y-1)^2:\ x,y\in\Z\}
\\\iff&n\in\{p_5(x)+2p_5(y):\ x,y\in\Z\},
\endalign$$
and
$$\align n\in\{T_x+2T_y:\ x,y\in\Z\}
\iff&8n+3\in\{u^2+2v^2:\ u,v\in\Z\}
\\\iff&24n+9\in\{x^2+2y^2:\ x,y\in\Z\ \&\ 3\nmid xy\}
\\\iff&24n+9\in\{x^2+8y^2:\ x,y\in\Z,\ 2\nmid x\ \&\ 3\nmid xy\}
\\\iff&24n+9\in\{(6x-1)^2+8(3y-1)^2:\ x,y\in\Z\}
\\\iff& n=p_5(x)+p_8(y)\ \t{for some}\ x,y\in\Z.
\endalign$$
This proves (1.7) and (1.8).

Now we show (1.10). Clearly,
$8(T_x+T_y)+2=(2x+1)^2+(2y+1)^2$
and
$$40(\psi_{5,1}(x)+\psi_{5,3}(y))+10=(10x+1)^2+(10y+3)^2=5((4x+2y+1)^2+(4y-2x+1)^2).$$
If $8n+2=u^2+v^2$ with $u$ and $v$ odd, then $40n+10=(2u+v)^2+(u-2v)^2=s^2+t^2$
for some $s,t\in\Z$ with $5\nmid st$ (by [S17, Lemma 2.1]), and hence $40n+10=(10x+1)^2+(10y+3)^2$ for some $x,y\in\Z$.
Therefore (1.10) holds.

The proof of Theorem 1.9 is now complete. \qed

\proclaim{Lemma 5.1} Let $w=3u^2+5v^2\in\Z^+$ with $u,v\in\Z$ and $8\mid w$. Then $w=3x^2+5y^2$ for some odd integers $x$ and $y$.
\endproclaim
\Proof. Let $k=\ord_2\gcd(u,v)$ and write $u=2^ku_0$ and $v=2^kv_0$ with $u_0,v_0\in\Z$ not all even.
If $k\in\{0,1\}$, then $u_0$ and $v_0$ are both odd since $8\mid w$. If $u_0\not\eq v_0\pmod 2$, then
$k\gs2$ and $4^2(3u_0^2+5v_0^2)=3u_2^2+5v_2^2$ with $u_2=u_0-5v_0$ and $v_2=3u_0+v_0$ both odd.

Let $j\in\N$. If $4^j(3u_0^2+5v_0^2)$ can be written as $3u_j^2+5v_j^2$ with $u_j$ and $v_j$ odd, then we may assume $u_j\not\eq v_j\pmod 4$
without loss of generality, hence
$$4^{j+1}(3u_0^2+5v_0^2)=4(3u_j^2+5v_j^2)=3u_{j+1}^2+5v_{j+1}^2$$
with $u_{j+1}=(v_j-u_j)/2+2v_j$ and $v_{j+1}=(v_j-u_j)/2+2u_j$ both odd.

By the above, $w=4^k(3u_0^2+5v_0^2)=3u_k^2+5v_k^2$ for some odd integers $u_k$ and $v_k$.  \qed
\medskip

\Remark\ 5.2. Note also the following useful fact:
$$3\l(\f x2+y\r)^2+5\l(\f x2-y\r)^2=3\l(\f{x-3y}2\r)^2+5\l(\f{x+y}2\r)^2=2x^2-2xy+8y^2.\tag5.1$$

The following lemma is a well known result in the theory of quadratic forms.

\proclaim{Lemma 5.3} {\rm ([C, Theorem 1.3])} Let $f$ be an integral quadratic form with nonzero discriminant. If an integer $m$ is represented by $f$ over the field of real numbers
as well as the ring $\Z_p$ of $p$-adic integers for each prime $p$, then $m$ is represented over $\Z$ by some form $f^*$ in the same genus as $f$.
\endproclaim

\proclaim{Lemma 5.4} Let $n\in\N$ and $\da\in\{0,1\}$. Then
$12n+8+3\da\in\{3x^2+3y^2+5z^2:\ x,y,z\in\Z\}.$
\endproclaim
\Proof. There are two classes in the genus of $3x^2+3y^2+5z^2$, the one not containing $3x^2+3y^2+5z^2$ has the representative
$3x^2+2y^2+8z^2-2yz$.
If $12n+8+3\da=3x^2+2y^2+8z^2-2yz$ with $x,y,z\in\Z$, $2\nmid y$ and $y\not\eq z\pmod2$, then $3\da\eq 3x^2+2y^2\eq -x^2+2\pmod 4$
which is impossible. Combining this with (5.1) and Lemma 5.3, we immediately obtain the desired result. \qed
\medskip

\noindent{\it Proof of Theorem 1.11}. Fix a nonnegative integer $n$.

(i) It is easy to see that
$$\align &n=\f{x(x+1)}2+\f{y(3y+1)}2+\f{z(5z+1)}2
\\\iff&120n+23=15(2x+1)^2+5(6y+1)^2+3(10z+1)^2.
\endalign$$

There are two classes in the genus of $3x^2+5y^2+15z^2$, and the one not containing $3x^2+5y^2+15z^2$ has the representative
$2x^2-2xy+8y^2+15z^2$.
If $120n+23=2x^2+8y^2+15z^2-2xy$ for some $x,y\in\Z$ with $2\nmid x$ and $y\not\eq x\pmod2$, then $23\eq 2x^2+15z^2\eq 17\pmod 4$ which is impossible. Thus, in view of (5.1) and Lemma 5.3, there are $x,y,z\in\Z$ such that $120n+23=3x^2+5y^2+15z^2$.

If $2\nmid x$, then $5(y^2+3z^2)\eq 23-3x^2\eq20\pmod 8$ and hence $y^2+3z^2=s^2+3t^2$ for some odd integers $s$ and $t$ (cf. [S15, Lemma 3.2]).
If $2\nmid z$, then $3x^2+5y^2\eq 23-15z^2\eq0\pmod 8$
and hence $3x^2+5y^2=3u^2+5v^2$ for some odd integers $u$ and $v$ (by Lemma 5.1).
If $x$ and $z$ are both even, then $y^2\eq 5y^2\eq23\eq3\pmod 4$ which is impossible.
So we may simply assume $2\nmid xyz$ without loss of generality.

Since $3x^2\eq 23\eq3\pmod 5$, $x$ or $-x$ is congruent to 1 modulo 10. As $y\not\eq0\pmod3$, $y$ or $-y$ is congruent to 1 modulo 6.
Thus, for some $u,v,w\in\Z$ we have
$$120n+23=3(10w+1)^2+5(6v+1)^2+15(2u+1)^2$$
and hence $ n=u(u+1)/2+v(3v+1)/2+w(5w+1)/2.$
This ends our proof of (1.12).

(ii) Let $\da\in\{0,1\}$ and $r\in\{1,3,5\}$.
There are two classes in the genus of $3x^2+5y^2+30z^2$, and the one not containing $3x^2+5y^2+30z^2$ is
$$2x^2+15y^2+15z^2=2x^2+30\l(\f{y+z}2\r)^2+30\l(\f{y-z}2\r)^2.$$
When $120n+30\da+3r^2+5=2x^2+30u^2+30v^2$ with $x,u,v\in\Z$, if $u\eq v\pmod2$ then $x\eq u\eq v\eq\da\pmod2$ (since $2\mid x-\da$ and $u^2+v^2\eq2\da\pmod4$),
thus we may assume $x\eq u\pmod2$ without loss of generality, and hence $2x^2+30u^2=3a^2+5b^2$
with $a=(x+5u)/2$ and $b=(x-3u)/2$ both integral.
So, with the help of Lemma 5.3, there are $x,y,z\in\Z$ such that $120n+30\da+3r^2+5=3x^2+5y^2+30z^2$.

Clearly, $z=2w+\da$ for some $w\in\Z$.  Since $3x^2+5y^2\eq 0\pmod 8$ and $3x^2+5y^2\not=0$, by Lemma 5.1 we can write
$3x^2+5y^2=3s^2+5t^2$ with $s$ and $t$ odd. Now, $120n+30\da+3r^2+5=3s^2+5t^2+120w(w+\da)+30\da$.
As $3s^2\eq 3r^2\pmod 5$, $s$ or $-s$ is congruent to $r$ modulo 10. Also, $t$ or $-t$ is congruent to $1$ modulo 6.
So there are $u,v\in\Z$ such that
$$120n+3r^2+5=3(10v+r)^2+5(6u+1)^2+120w(w+\da)$$
and hence
$n=w(w+\da)+u(3u+1)/2+v(5v+r)/2.$
This proves (1.13).

(iii) Let $r,s,t\in\{1,3\}$ with $\{r,s\}\not=\{3\}$.
There are two classes in the genus of $3x^2+5y^2+5z^2$, and the one not containing $3x^2+5y^2+5z^2$ has the representative
$2x^2-2xy+8y^2+5z^2$.
If $120n+5r^2+5s^2+3t^2=2x^2+8y^2-2xy+5z^2$ with $x,y,z\in\Z$, $2\nmid x$ and $2\mid y$, then $13\eq 5r^2+5s^2+3t^2\eq 2x^2+5z^2\eq 2+5\pmod 4$
which is impossible. So, in light of (5.1) and Lemma 5.3, there are $x,y,z\in\Z$ such that $120n+5r^2+5s^2+3t^2=3x^2+5y^2+5z^2$. As $3x^2\not\eq 13\eq 5r^2+5s^2+3t^2\pmod4$,
$y$ and $z$ cannot be both even. Without loss of generality, we assume that $2\nmid z$. Then $3x^2+5y^2>0$ and $3x^2+5y^2\eq0\pmod 8$.
By Lemma 5.1, we can write $3x^2+5y^2$ as $3x_0^2+5y_0^2$ with $x_0$ and $y_0$ both odd.

By the last paragraph,  $120n+5r^2+5s^2+3t^2=3x^2+5y^2+5z^2$ for some odd integers $x,y,z$.
Clearly $x$ or $-x$ has the form $10w+t$ with $w\in\Z$. Since $y^2+z^2\eq r^2+s^2\pmod3$,
we have
$y^2+z^2=(6u+r)^2+(6v+s)^2$ for some $u,v\in\Z$. Therefore
$$120n+5r^2+5s^2+3t^2=3(10w+t)^2+5(6u+r)^2+5(6v+s)^2$$
and hence
$n=u(3u+r)/2+v(3v+s)/2+w(5w+t)/2.$
This proves (1.14).

(iv) As $3(s^2+t^2)+5\eq 3\times2+5=11\pmod{12}$, by Lemma 5.4 there are $x,y,z\in\Z$ such that $120n+3(s^2+t^2)+5=3x^2+3y^2+5z^2$.
Clearly, $x$ and $y$ cannot be both even. Without loss of generality, we assume that $2\nmid x$. Then $3y^2+5z^2>0$ and
$3y^2+5z^2\eq 0\pmod 8$.
By Lemma 5.1, we can write $3y^2+5z^2$ as $3y_0^2+5z_0^2$ with $y_0$ and $z_0$ both odd.
So, without loss of generality we may simply assume that $y$ and $z$ are also odd.

If $\{s,t\}=\{1,3\}$, then $\psi_{5,s}(y)+\psi_{5,t}(z)\sim T_y+T_z$ by (1.10), hence
$\psi_{3,1}(x)+\psi_{5,s}(y)+\psi_{5,t}(z)$ is universal over $\Z$ as $p_5(x)+T_y+T_z$ is universal over $\Z$ by [S15, Theorem 1.14].

Now we assume that $\{s,t\}\not=\{1,3\}$.
Clearly, $z$ or $-z$ has the form $6w+1$ with $w\in\Z$. Since $x^2+y^2\eq s^2+t^2\not\eq0\pmod5$, we have
$x^2+y^2=(10u+s)^2+(10v+t)^2$ for some $u,v\in\Z$. Therefore
$$120n+3s^2+3t^2+5=5(6w+1)^2+3(10u+s)^2+3(10v+t)^2$$
and hence $n=w(3w+1)/2+u(5u+s)/2+v(5v+t)/2.$
This proves (1.15).

In view of the above, we have completed the proof of Theorem 1.11. \qed

\medskip

\Ack. The author would like to thank the two referees and his graduate student Hai-Liang Wu for helpful comments.

\newpage
\heading{Appendix} \endheading

In this appendix, we list our conjectural universal tuples $(a,b,c,d,e,f)$ over $\N$ with $a\gs c\gs e\gs1$, $b>-a$ and $b\eq a\pmod2$,
$d>-c$ and $d\eq c\pmod2$, $f>-e$ and $f\eq e\pmod2$, and ($a\nmid b$ or $c\nmid d$ or $e\nmid f$).
For each of the listed tuple $(a,b,c,d,e,f)$, we have verified that every $n=0,\ldots,10^6$ can be written as $\psi_{a,b}(x)+\psi_{c,d}(y)+\psi_{e,f}(z)$
with $x,y,z\in\N$. Below is our list.
$$\align &(3, -1, 1, 1, 1, 1),\ (3,-1, 1, 3, 1, 1),\ (3, -1, 1, 5, 1, 1),\ (3, -1, 1, 7, 1, 1),
\\& (3, -1, 1, 9, 1, 1),\ (3, -1, 1, 11, 1, 1),\ (3, -1, 1, 13, 1, 1),\ (3, -1, 1, 15, 1, 1),
\\&(3, -1, 1, 17, 1, 1),\ (3, -1, 2, 0, 1, 1),\ (3, -1, 2, 0, 1, 3),\ (3, -1, 2, 0, 1, 5),
\\&(3, -1, 2, 0, 2, 0),\ (3, -1, 2, 2, 1, 1),\ (3, -1, 2, 2, 1, 3),\ (3, -1, 2, 2, 2, 0),
\\&(3, -1, 2, 4, 1, 1),\ (3, -1, 2, 4, 1, 3),\ (3, -1, 2, 4, 2, 2),\ (3, -1, 2, 6, 1, 1),
\\&(3, -1, 2, 8, 1, 1),\ (3, -1, 2, 10, 1, 1),\ (3, -1, 2, 12, 1, 1),\ (3, -1, 2, 16, 1, 1),
\\&(3, -1, 3, -1, 1, 1),\ (3, -1, 3, -1, 1, 3),\ (3, 1, 1, 1, 1, 1),\ (3, 1, 1, 3, 1, 1),
\\&(3, 1, 1, 5, 1, 1),\ (3, 1, 1, 7, 1, 1),\ (3, 1, 2, 0, 1, 1),\ \ (3, 1, 2, 0, 1, 3),
\\&(3, 1, 2, 2, 1, 1),\ (3, 1, 2, 2, 2, 0),\ (3, 1, 2, 4, 1, 1),\ (3, 1, 2, 6, 1, 1),
\\&(3, 1, 2, 6, 2, 0),\ (3, 1, 3, -1, 1, 1),\ (3, 1, 3, -1, 1, 3),\ (3, 1, 3, -1, 2, 2),
\\&(3, 1, 3, -1, 2, 6),\ (3, 3, 3, -1, 1, 1), (3, 3, 3, -1, 1, 3),\ (3, 3, 3, -1, 2, 0),
\\&(3, 3, 3, -1, 2, 2),\ (3, 3, 3, 1, 1, 1),\ (3, 3, 3, 1, 2, 0),\ (3, 5, 1, 1, 1, 1),
\\&(3, 5, 1, 3, 1, 1),\ (3, 5, 2, 0, 1, 1),\ (3, 5, 2, 0, 1, 3),\ (3, 5, 2, 2, 1, 1),
\\&(3, 5, 2, 2, 2, 0),\ (3, 5, 3, -1, 1, 1),\ (3, 5, 3, 1, 1, 1),\ (3, 7, 1, 1, 1, 1),
\\&(3, 7, 2, 0, 1, 1),\ (3, 7, 2, 2, 2, 0),\ (3, 7, 3, -1, 1, 1),\ (3, 9, 3, -1, 1, 1),
\\&(3, 11, 2, 0, 1, 1),\ (3, 11, 2, 0, 1, 3),\ (3, 11, 3, -1, 1, 1),\ (3, 13, 2, 0, 1, 1),
\\&(3, 15, 3, -1, 1, 1);
\endalign$$
$$\align & (4, -2, 1, 1, 1, 1),\ (4, -2, 1, 3, 1, 1),\ (4, -2, 1, 3, 1, 3),\ (4, -2, 1, 5, 1, 1),
\\& (4, -2, 1, 5, 1, 3),\ (4, -2, 1, 7, 1, 1),\ (4, -2, 1, 7, 1, 3),\ (4, -2, 1, 9, 1, 1),
\\& (4, -2, 2, 0, 1, 1),\ (4, -2, 2, 0, 1, 3),\ (4, -2, 2, 0, 1, 5),\ (4, -2, 2, 2, 1, 1),
\\& (4, -2, 2, 2, 1, 3),\ (4, -2, 2, 2, 1, 5),\ (4, -2, 2, 2, 2, 0),\ (4, -2, 2, 4, 1, 1),
\\& (4, -2, 2, 8, 1, 1),\ (4, -2, 3, -1, 1, 1),\ \ (4, -2, 3, -1, 1, 3),\ (4, -2, 3, -1, 1, 5),
\\& (4, -2, 3, -1, 2, 0),\ (4, -2, 3, 1, 1, 1), \ (4, -2, 3, 1, 1, 3),\ (4, -2, 3, 1, 2, 0),
\\&(4, -2, 3, 7, 1, 1),\ (4, 0, 3, -1, 1, 1),\ (4, 0, 3, -1, 1, 3),\ (4, 0, 3, -1, 1, 5),
\\& (4, 0, 3, -1, 1, 7),\ (4, 0, 3, -1, 2, 0),\ (4, 0, 3, -1, 2, 4),\ (4, 0, 3, 1, 1, 1),
\\&(4, 0, 3, 1, 2, 0),\ (4, 0, 3, 5, 1, 1),\ (4, 0, 3, 5, 2, 0),\ (4, 0, 4, -2, 1, 1).
\\&(4, 0, 4, -2, 1, 3),\ (4, 0, 4, -2, 1, 5),\ (4, 0, 4, -2, 3, -1),\ (4, 2, 2, 0, 1, 1).
\\&(4, 2, 2, 0, 1, 3),\ (4, 2, 3, -1, 1, 1),\ (4, 2, 3, -1, 1, 3),\ (4, 4, 3, -1, 1, 1),
\\&(4, 4, 3, 1, 1, 1),\ (4, 6, 3, -1, 1, 1),\ (4, 8, 3, -1, 1, 1),\ (4, 10, 3, -1, 1, 1);
\endalign$$
$$\align
\\&(5, -3, 1, 1, 1, 1),\ (5, -3, 1, 3, 1, 1),\ (5, -3, 1, 5, 1, 1),\ (5, -3, 1, 7, 1, 1),
\\&(5, -3, 1, 9, 1, 1),\ (5, -3, 2, 0, 1, 1),\ (5, -3, 2, 0, 1, 3),\ (5, -3, 2, 2, 1, 1),
\\& (5, -3, 2, 2, 1, 3),\ (5, -3, 2, 2, 2, 0),\ (5, -3, 2, 4, 1, 1),\ (5, -3, 2, 4, 1, 3),
\\&(5, -3, 2, 4, 2, 2),\ (5, -3, 2, 8, 1, 1),\ (5, -3, 3, -1, 1, 1),\ (5, -3, 3, 1, 1, 1),
\\&(5, -3, 3, 1, 1, 3),\ (5, -3, 3, 3, 3, -1),\ (5, -3, 3, 7, 1, 1),\ (5, -1, 2, 0, 1, 1),
\\&(5, -1, 2, 0, 1, 3),\ (5, -1, 2, 0, 1, 5),\ (5, -1, 2, 0, 1, 7),\ (5, -1, 2, 0, 1, 9),
\\& (5, -1, 2, 2, 1, 1),\ (5, -1, 2, 6, 1, 1),\ (5, -1, 3, -1, 1, 1),\ (5, -1, 3, 1, 1, 1),
\\&(5, -1, 4, 0, 1, 1),\ (5, 1, 2, 0, 1, 1),\ \ (5, 1, 2, 0, 1, 3),\ (5, 1, 2, 2, 1, 1),
\\&(5, 1, 3, -1, 1,1),\ (5, 1, 3, -1, 1, 3),\ (5, 1, 3, 1, 1, 1),\ (5, 3, 1, 1, 1, 1),
\\&(5, 3, 1, 3, 1, 1),\ (5, 3, 2, 0, 1, 1),\ (5, 3, 2, 2, 1, 1),\ (5, 3, 2, 2, 2, 0),
\\&(5, 3, 3, -1, 1, 1),\ (5, 3, 3, 1, 2, 0),\ (5, 3, 4, -2, 1, 1),\ (5, 3, 4, -2, 1, 3),
\\&(5, 7, 2, 0, 1, 1),\ (5, 7, 2, 0, 1, 3),\ (5, 7, 3, -1, 1, 1),\ (5, 9, 3, -1, 1, 1),
\\&(5, 11, 2, 0, 1, 1),\ (5, 13, 3, -1, 1, 1);
\endalign$$
$$\align&(6, -4, 1, 1, 1, 1),\ (6, -4, 1, 3, 1, 1),\ (6, -4, 1, 3, 1, 3),\ (6, -4, 1, 5, 1, 1),
\\&(6, -4, 1, 5, 1, 3),\ (6, -4, 1, 7, 1, 1),\ (6, -4, 1, 7, 1, 3),\ (6, -4, 1, 9, 1, 1),
\\&(6, -4, 2, 0, 1, 1),\ (6, -4, 2, 0, 1, 3),\ (6, -4, 2, 2, 1, 1),(6, -4, 2, 2, 1, 3),
\\&(6, -4, 2, 2, 2, 0),\ \ (6, -4, 2, 6, 1, 1),\ (6, -4, 2, 6, 1, 3),\ (6, -4, 3, -1, 1,1),
\\&(6, -4, 3, 1, 1, 1),\ (6, -4, 3, 1, 1, 3),\ (6, -4, 3, 1, 1, 5),\ (6, -4, 3, 5, 1, 1),
\\&(6, -4, 3, 7, 1, 1),\ (6, -4, 4, 0, 1, 1),\ (6, -4, 4, 0, 1, 3),\ (6, -4, 4, 0, 1, 5),
\\&(6, -4, 4, 0, 1, 7),\ (6, -4, 5, -3, 1, 1),\ (6, -2, 1, 1, 1, 1),\ (6, -2, 1, 3, 1, 1),
\\&(6, -2, 1, 5, 1, 1),\ (6, -2, 1, 7, 1, 1),\ (6, -2, 2, 0, 1, 1),\ (6, -2, 2, 0, 1, 3),
\\&(6, -2, 2, 0, 1, 5),\ (6, -2, 2, 0, 1, 7),\ (6, -2, 2, 0, 1, 9),\ (6, -2, 2, 2, 2, 0),
\\&(6, -2, 3, 1, 1, 1),\ (6, -2, 4, -2, 1, 1),\ (6, -2, 4, -2, 1, 3),\ (6, -2, 4, -2, 1, 5),
\\&(6, -2, 4, 0, 1, 1),\ (6, -2, 4, 4, 1, 1),\ (6, -2, 5, -3, 1, 1),\ (6, -2, 5, -3, 1, 3),
\\&(6, 0, 3, -1, 1, 1),\ (6, 0, 3, -1, 2, 2),\ (6, 0, 3, 1, 1, 1),\ (6, 0, 6, -4, 1, 1),
\\&(6, 0, 6, -4, 2, 2),\ (6, 2, 1, 1, 1, 1),\ (6, 2, 1, 3, 1, 1),\ (6, 2, 2, 2, 1, 1),
\\&(6, 2, 2, 2, 2, 0),\ (6, 2, 3, -1, 1, 1),\ (6, 2, 3, 1, 1, 1),\ (6, 2, 4, 0, 1, 1),
\\&(6, 2, 6, -4, 1, 1),\ (6, 4, 1, 1, 1, 1),\ (6, 4, 2, 2, 2, 0),\ (6, 4, 6, -4, 1, 1),
\\&(6, 6, 3, -1, 1, 1),\ (6, 8, 2, 0, 1, 1),\ (6, 8, 3, -1, 1, 1),\ (6, 10, 2, 0, 1, 1),
\\&(6, 12, 3, -1, 1, 1);
\endalign$$
$$\align&(7, -5, 2, 0, 1, 1),\ (7, -5, 2, 2, 1, 1),\ (7, -5, 2, 4, 1, 1),\ (7, -5, 2, 6, 1, 1),
\\& (7, -5, 2, 8, 1, 1),\ (7, -5, 3, -1, 1, 1),\ (7, -5, 3, -1, 1, 3),\ (7, -5, 3, -1, 2, 2),
\\& (7, -5, 3, 5, 1, 1),\ (7, -5, 4, 0, 1, 1),\ (7, -3, 2, 0, 1, 1),\ (7, -3, 2, 2, 1, 1),
\\& (7, -3, 2, 6, 1, 1),\ (7, -3, 3, -1, 1, 1),\ (7, -3, 3, -1, 1, 3),\ (7, -3, 3, -1, 1, 5),
\\&(7, -3, 3, 1, 1, 1),\ (7, -3, 4, 0, 1, 1),\ (7, -1, 2, 0, 1, 1),\ (7, -1, 2, 0, 1, 3),
\\&(7, -1, 2, 2, 1, 1),\ (7, -1, 3, 1, 1, 1),\ (7, -1, 6, -4, 1, 1),\ (7, 1, 2, 0, 1, 1),
\\&(7, 1, 3, -1, 1, 1),\ (7, 3, 1, 1, 1, 1),\ (7, 3, 2, 2, 2, 0),\ (7, 5, 2, 0, 1, 1),
\\& (7, 5, 2, 0, 1, 3);
\endalign$$
$$\align& (8, -6, 1, 1, 1, 1),\ (8, -6, 1, 3, 1, 1),\ (8, -6, 1, 5, 1, 1),\ (8, -6, 1, 7, 1, 1),
\\&(8, -6, 2, 0, 1, 1),\ (8, -6, 2, 2, 1, 1),\ (8, -6, 2, 2, 1, 3),\ (8, -6, 2, 2, 2, 0),
\\&(8, -6, 2, 4, 2, 2),\ (8, -6, 3, 1, 1, 1),\ (8, -6, 5, -3, 1, 1),\ (8, -6, 5, -3, 1, 3),
\\& (8, -4, 2, 2, 1, 1),\ (8, -4, 2, 6, 1, 1),\ (8, -4, 6, -4, 1, 1),\ (8, -2, 3, -1, 1, 1),
\\&(8, -2, 3, 1, 1, 1),\ (8, 0, 3, -1, 1, 1),\ (8, 0, 3, 1, 1, 1),\ (8, 0, 6, -2, 1, 1),
\\&(8, 2, 1, 1, 1, 1),\ (8, 2, 2, 2, 2, 0);
\endalign$$
$$\align&(9, -7, 2, 0, 1, 1),\ (9, -7, 2, 4, 1, 1),\ (9, -7, 2, 8, 1, 1),\ (9, -7, 3, -1, 1, 1),
\\&(9, -7, 5, -1, 1, 1),\ (9, -7, 5, -1, 1, 3),\ (9, -5, 2, 2, 1, 1),\ (9, -5, 2, 6, 1, 1),
\\&(9, -5, 3, -1, 1, 1),\ (9, -5, 3, -1, 1, 3),\ (9, -5, 3, 5, 1, 1),\ (9, -5, 4, 0, 1, 1),
\\&(9, -1, 1, 1, 1, 1),\ (9, -1, 1, 3, 1, 1),\ (9, -1, 2, 0, 1, 1),\ (9, -1, 2, 2, 1, 1),
\\&(9, -1, 2, 2, 2, 0),\ (9, -1, 4, 0, 1, 1),\ (9, 5, 3, -1, 1, 1),\ (9, 7, 2, 0, 1, 1),
\\&(9, 9, 3, -1, 1, 1);
\endalign$$
$$\align&(10, -8, 1, 1, 1, 1),\ (10, -8, 1, 3, 1, 1),\ (10, -8, 1, 5, 1, 1),\ (10, -8, 1, 7, 1, 1),
\\&(10, -8, 1, 9, 1, 1),\ (10, -8, 2, 0, 1, 1),\ (10, -8, 2, 2, 1, 1),\ (10, -8, 2, 2, 1, 3),
\\& (10, -8, 2, 2, 2, 0),\ (10, -8, 2, 4, 1, 1),\ (10, -8, 2, 8, 1, 1),\ (10, -8, 3, 1, 1, 1),
\\& (10, -8, 3, 1, 1, 3),\ \ (10, -8, 3, 7, 1, 1),\ (10, -6, 2, 0, 1, 1),\ (10, -6, 2, 0, 1, 3),
\\&(10, -6, 2, 0, 1, 5),\ (10, -6, 2, 2, 1, 1),\ (10, -6, 3, -1, 1, 1),\ (10, -6, 3, 5, 1, 1),
\\&(10, -6, 5, -3, 1, 1),\ (10, -4, 2, 0, 1, 1),\ (10, -4, 2, 0, 1, 3),\ (10, -4, 3, -1, 1, 1),
\\& (10, -2, 3, -1, 1, 1),\ (10, 2, 3, -1, 1, 1),\ (10, 4, 3, -1, 1, 1),\ (10, 6, 2, 0, 1, 1),
\\& (10, 8, 3, -1, 1, 1);
\endalign$$
$$\align
\\&(11, -9, 2, 0, 1, 1),\ (11, -9, 3, -1, 1, 1),\ (11, -9, 3, 1, 1, 1),\ (11, -9, 3, 7, 1, 1),
\\&(11, -7, 1, 1, 1, 1),\ (11, -7, 1, 3, 1, 1),\ (11, -7, 1, 5, 1, 1),\ (11, -7, 2, 2, 2, 0),
\\& (11, -7, 3, -1, 1, 1),\ (11, -5, 2, 0, 1, 1),\ (11, -5, 2, 0, 1, 3),\ (11, -3, 3, -1, 1, 1),
\\&(11, -3, 3, 1, 1, 1),\ (11, 5, 2, 0, 1, 1);
\endalign$$
$$\align&(12, -8, 2, 0, 1, 1),\ (12, -8, 2, 0, 1, 3),\ (12, -8, 2, 2, 1, 1),\ (12, -4, 1, 1, 1, 1),
\\&(12, -4, 1, 3, 1, 1),\ (12, -4, 2, 2, 1, 1),\ (12, -4, 2, 2, 2, 0),\ (12, -4, 4, 0, 1, 1).
\\& (13, -11, 2, 0, 1, 1),\ (13, -11, 2, 2, 1, 1),\ (13, -11, 2, 2, 1, 3),\ (13, -9, 2, 0, 1, 1),
\\& (13, -9, 2, 0, 1, 3), \ (13, -7, 2, 0, 1, 1),\ (13, -7, 2, 0, 1, 3),\ (13, -5, 2, 0, 1, 1).
\\&(14, -12, 2, 2, 1, 1),\ (14, -10, 2, 0, 1, 1),\ (14, -10, 2, 0, 1, 3),\ (14, -10, 2, 2, 1, 1),
\\&(14, -10, 3, -1, 1, 1),\ (14, -2, 2, 0, 1, 1),\ (14, -2, 3, -1, 1, 1),\ (14, 2, 2, 0, 1, 1),
\\&(15, -13, 1, 1, 1, 1),\ (15, -13, 1, 3, 1, 1),\ (15, -13, 1, 5, 1, 1),\ (15, -13, 1, 7, 1, 1),
\\&(15, -13, 1, 9, 1, 1),\ (15, -13, 2, 0, 1, 1),\ (15, -13, 2, 2, 1, 1),\ (15, -13, 2, 2, 1, 3),
\\&(15, -13, 2, 2, 2, 0),\ (15, -13, 2, 4, 1, 1),\ (15, -13, 2, 4, 2, 2),\ (15, -13, 2, 8, 1, 1),
\\&(15, -13, 3, 1, 1, 1),\ (15, -13, 3, 7, 1, 1),\ (15, -13, 4, 0, 1, 1),\ (15, -11, 1, 1, 1, 1),
\\&(15, -11, 2, 0, 1, 1),\ (15, -11, 2, 2, 2, 0),\ (15, -7, 1, 1, 1, 1),\ (15, -7, 2, 0, 1, 1),
\\&(15, -7, 2, 2, 1, 1),\ (15, -7, 2, 2, 2, 0),\ (15, -3, 3, -1, 1, 1),\ (15, 3, 3, -1, 1, 1),
\\& (16, -14, 2, 0, 1, 1),\ (16, -10, 2, 0, 1, 1),\ (16, -10, 2, 0, 1, 3),\ (16, -8, 3, -1, 1, 1),
\\& (16, -4, 2, 0, 1, 1),\ (17, -15, 3, 1, 1, 1),\ (17, -15, 3, 1, 1, 3),\ (18, -10, 2, 0, 1, 1),
\\& (20, -16, 2, 2, 1, 1),\ (20, -16, 2, 6, 1, 1),\ (20, -12, 3, -1, 1, 1),\ (20, -4, 2, 0, 1, 1),
\\&(21, -19, 2, 2, 1, 1),\ (21, -9, 2, 0, 1, 1),\ (21, -5, 2, 0, 1, 1),\ (25, -23, 2, 0, 1, 1).
\endalign$$

\enddocument